\documentclass{amsart}
\usepackage{amssymb}
\usepackage{verbatim}
\usepackage{array}
\usepackage{latexsym}
\usepackage{enumerate}
\usepackage{amsmath}
\usepackage{amsfonts}
\usepackage{amsthm}
\usepackage{color}
\usepackage[english]{babel}
\usepackage{comment}
\usepackage[pagewise]{lineno}
\newtheorem{theorem}{Theorem}[section]
\newtheorem{proposition}[theorem]{Proposition}
\newtheorem{lemma}[theorem]{Lemma}

\newtheorem{result}[theorem]{Result}

\def\cC{\mathcal C}
\def\cD{\mathcal D}

\def\cF{\mathcal F}

\def\cH{\mathcal H}
\def\cK{\mathcal K}

\def\cQ{\mathcal Q}

\def\cS{\mathcal S}
\def\cU{\mathcal U}

\def\cW{\mathcal W}

\def\PG{{\rm{PG}}}

\def\deg{\mbox{\rm deg}}

\def\Sym{\mbox{\rm Sym}}

\newcommand{\PGL}{\mbox{\rm PGL}}
\newcommand{\AGL}{\mbox{\rm AGL}}

\newcommand{\Gal}{{\rm {Gal}}}
\newcommand{\ha}{{\textstyle\frac{1}{2}}}

\title[Algebraic approach to the completeness problem for $(k,n)$-arcs]{Algebraic approach to the completeness problem for $(k,n)$-arcs in planes over finite fields}
\date{}

\author{G\'abor Korchm\'aros}
\address{Dipartimento di Matematica, Informatica ed Economia\\
	Universit\`a degli Studi della Basilicata\\
	viale dell'Ateneo Lucano 10\\
	85100 Potenza, Italy}
\email{gabor.korchmaros@unibas.it}
\author{G\'abor P. Nagy}
\address{Bolyai Institute \\
	University of Szeged \\
	Aradi v\'ertan\'uk tere 1\\
	H-6720 Szeged, Hungary}
\address{Department of Algebra \\
	Budapest University of Technology and Economics\\
	M\H{u}egyetem rkp. 3\\
	H-1111 Budapest, Hungary}
\email{nagyg@math.u-szeged.hu, nagy.gabor.peter@ttk.bme.hu}
\author{Tam\'as Sz\H onyi}
\address{Institute of Mathematics, ELTE E\"otv\"os Lor\'and University and MTA-ELTE Geometric and Algebraic Combinatorics Research Group \\
	P\'azm\'any P\'eter s\'et\'any 1/C \\
	H-1117 Budapest, Hungary}
\address{FAMNIT, University of Primorska \\
	Glagoljaška 8 \\
	6000 Koper, Slovenia}
\email{tamas.szonyi@ttk.elte.hu, tamas.szonyi@famnit.upr.si}

\begin{document}

%\thanks{{\em Keywords}: Algebraic curves, Positive characteristic, Automorphism groups.}

\begin{abstract} In a projective plane over a finite field, complete $(k,n)$-arcs with few characters are rare but interesting objects with several applications to finite geometry and coding theory. Since almost all known examples are large, the construction of small ones, with $k$ close to the order of the plane, is considered a hard problem. A natural candidate to be a small $(k,n)$-arc with few characters is the set $\Omega(\cC)$ of the points of a plane curve $\cC$ of degree $n$ (containing no linear components) such that some line meets $\cC$ transversally in the plane, i.e. in $n$ pairwise distinct points. Let $\cC$ be either the Hermitian curve of degree $q+1$ in $\PG(2,q^{2r})$ with $r\ge 1$, or the rational BKS curve of degree $q+1$ in $\PG(2,q^r)$ with $q$ odd and $r\ge 1$. Then $\Omega(\cC)$ has four and seven characters, respectively. Furthermore, $\Omega(\cC)$ is small as both curves are either maximal or minimal. The completeness problem is investigated by an algebraic approach based on Galois theory and on the Hasse-Weil lower bound. Our main result for the Hermitian case is that $\Omega(\cC)$ is complete for $r\ge 4$. For the rational BKS curve, $\Omega(\cC)$  is complete if and only if $r$ is even. If $r$ is odd then the uncovered points by the $(q+1)$-secants to $\Omega(\cC)$ are exactly the points in $\PG(2,q)$ not lying in $\Omega(\cC)$. Adding those points to $\Omega(\cC)$ produces a complete $(k,q+1)$-arc in $\PG(2,q^r)$, with $k=q^r+q$. The above results do not hold true for $r=2$ and there remain open the case $r=3$ for the Hermitian curve, and the cases $r=3,4$ for the rational  BKS curve. As a by product we also obtain two results of interest in the study of the Galois inverse problem for $\PGL(2,q)$.
\end{abstract}
\maketitle
\vspace{0.5cm}\noindent {\em Keywords}:
$(k,n)$-arcs in $\PG(2,q)$, algebraic curves, Galois theory
\vspace{0.2cm}\noindent

\vspace{0.5cm}\noindent {\em Subject classifications}:
\vspace{0.2cm}\noindent  51E21, 14G15,
11R32

\section{Introduction}
Let $k,n$ with $2\le n<k$ be two positive integers. In the projective plane $\PG(2,q)$ over a finite field $\mathbb{F}_q$ of order $q$, a $(k,n)$-arc is point-set $\cK$ of size $k$ such that $n$ is the maximum number of collinear points in $\cK$. The concept of a $(k,n)$-arc has a useful interpretation in coding theory as the incidence matrix of a $(k,n)$-arc of $\PG(2,q)$ is the parity check matrix of a $[k,3,k-n]_q$ almost MDS code which is non-extendible if and only if the corresponding $(k,n)$-arc is complete. Here, completeness means that the $(k,n)$-arc is maximal (not contained properly in a larger $(k',n)$-arc), in other words  each point off the $(k,n)$-arc in $\PG(2,q)$ is incident to some $n$-secant i.e. lines meeting the $(k,n)$-arc in exactly $n$ points. An (incidence) character $n_i$ of a $(k,n)$-arc $\cK$ is any integer $0\le n_i\le n$ such that some line meets $\cK$ in exactly $n_i$ points.  The foundation of the theory of $(k,n)$-arcs with special attention on their characters was laid down in the 1950s by B. Segre and A. Barlotti.
%; later on many contributions were given especially by the Italian, Belgian and Hungarian schools of finite geometry.

In the smallest case $n=2$, plenty of results and constructions are known, especially for larger complete $(k,2)$-arcs; see \cite{H}. In particular, the maximum size of a $(k,2)$-arc is $q + 1$ or $q+2$ according as $q$ odd or even. In the odd order case, the $(q+1,2)$-arcs are exactly the sets consisting of all points of an irreducible conic, whereas the classification of $(q+2,2)$-arcs is still open. The second and third largest $(k,2)$-arcs have also been studied intensively.
%On the other end, probabilistic methods show the existence of complete $(k,2)$-arc for $k=O(\sqrt{q}\log^c q)$ with a positive constant $c$. Several effective constructions for small $(k,2)$-arcs are available in the literature, the best explicit example so far is for $k=O(q^{3/4})$; see \cite{SzT}.
On the other end, there is an elementary lower bound for the size of a complete arc, due to Lunelli and Sce,  see \cite{H}, Thm. 9.12 and a slight improvement in Thm. 9.13. Probabilistic methods show the existence of complete $(k,2)$-arcs for $k=O(\sqrt q\log ^c q)$ with a positive constant $c$. Several effective constructions for small $(k,2)$-arcs are available in the literature (see e.g. \cite{SzT}), the best explicit example so far is for $k=O(q^{3/4})$.

%Our current knowledge of $(k,n)$-arcs with $n\ge 3$ is much less although well known combinatorial structures embedded in $\PG(2,q)$, such as unitals, nets and other designs, provide inspiring examples of $(k,n)$-arcs with few characters.
Our current knowledge of $(k,n)$-arcs with $n\ge 3$ is much less although well known combinatorial structures embedded in $\PG (2,q)$, such as unitals, nets and other designs, provide inspiring examples of $(k,n)$-arcs with few characters.The upper bound on $k$ for $(k,n)$-arcs is due to Barlotti (see \cite{H}, Corollary 12.5), and Ball, Blokhuis and Mazzocca \cite{BBM} showed that it cannot be achieved in $\PG (2,q)$, $q$ odd. More details can be found in Chapter 12 of \cite{H}. There is also a combinatorial lower bound, due to Alabdullah and Hirschfeld \cite{AH}, analogous to the Lunelli-Sce bound. Also in this case one can apply probabilistic methods, and also effective constructions, in particular when $n$ is reasonably large compared to $q$.
Since $(k,n)$-arcs are truly combinatorial objects, counting arguments and incidence geometry are prevalent in their studies. Algebraic and geometric methods combined with combinatorics have also been successfully developed to construct and investigate new and interesting families of $(k,n)$-arcs, the principal tools being polarities, geometric transformations, groups of symmetries, and algebraic curves. Nevertheless, the study of $(k,n)$-arcs is still hard whenever it requires that the $(k,n)$-arc be complete with only few characters and small size $k\approx q$.

A natural candidate to be such a complete $(k,n)$-arc is the set of the points of a plane curve $\cC$ of degree $n$ (containing no linear components over $\mathbb{F}_q$) such that some line of $\PG(2,q)$ meets $\cC$ transversally in $\PG(2,q)$, i.e. in $n$ pairwise distinct points in $\PG(2,q)$. The most interesting families of complete $(k,n)$-arcs with very few characters are of this kind including Baer subplanes, ovals, classical unitals and Denniston type maximal arcs; see Section \ref{bfg}. In particular, the classical unital in $\PG(2,q^2)$ consists of all points of the Hermitian curve $\cH_q$  which is also well known in number theory as being the most important $\mathbb{F}_{q^2}$-maximal curve. Further examples of $(k,n)$-arcs consisting of the points of (Frobenius non-classical) curves are found in \cite{GPUT} and \cite{BH,BMT}. The work by J.W.P. Hirschfeld and F.J. Voloch for cubic curves, and by M. Giulietti and F. Torres for $n\ge 4$, was the first important step towards a kind of sophisticated algebraic theory of $(k,n)$-arcs arising from plane algebraic curves.
A main result in this theory concerns plane curves of degree $n$ defined over a subfield $\mathbb{F}_{\bar{q}}$ of $\mathbb{F}_q$ and viewed as a curve of $\PG(2,q)$: If $q$ is large enough compared to the parameters of $\cC$, namely $n$ and $\bar{q}$, then the set of the points of any absolutely irreducible curve $\cC$ of degree $n$ in $\PG(2,q)$ is a $(k,n)$-arc of small size $k\approx q$.

The algebraic theory approach is also adequate for the completeness problem for such $(k,n)$-arcs even if it needs Galois theory in positive characteristic together with some Dirichlet or $\rm{\check{C}}$ebotarev type density theorem; see paper \cite{BM} that finds its origin in previous work by Guralnick, Tucker, Zieve \cite{GTZ} and others on permutation polynomials.  The essential idea is to express the condition that a point $P\in \PG(2,q)$ is incident with a line intersecting transversally $\cC$ in $\PG(2,q)$ in terms of the Galois closure of the algebraic extension $F|F_P$ where $F$ is the function field of $\cC$ and $F_P$ is the rational subfield of $F$ arising from the projection of $\cC$ from $P$. The favorable situation occurs when the (geometric) Galois group $Gal(F|F_P)$ is the symmetric group $\rm{\Sym}_n$ on the roots of the polynomial associated with $F|F_P$. In fact, for this case, a variant of the classical $\rm{\check{C}}$ebotarev density theorem \cite[Theorem 9.13B]{MR} works well and ensures the existence of a line $\ell$ through $P$ meeting transversally $\cC$ in $\PG(2,q)$ provided that $q$ is large enough compared to the two parameters of $\cC$, namely the degree of $\cC$ and the order of the plane $\PG(2,\bar{q})$  where $\cC$ is defined. For $Gal(F|F_P)\lvertneqq \rm{Sym}_n$, the $(k,n)$-arc may not be complete; nevertheless completeness can still be achieved by adding some (at most $O(n)$) points; see \cite{BM}. As a corollary, see \cite[Theorem 5.3]{BM}, for all but finitely many $n$'s, if $q$ is large enough, there are complete $(k,n)$-arcs of small size $k\approx q$.

The question arises whether complete $(k,n)$-arcs of small size $k\approx q$ can be obtained in this way for (almost) every $q$. As it appears plausible, at least intuitively, the choice of the curve is critical. We thoroughly work out two cases investigating the Hermitian curve $\cH_q$ defined over $\mathbb{F}_{q^2}$ and the rational BKS curve $\Gamma_q$ defined over $\mathbb{F}_q$, respectively. Our main result for the Hermitian case is that for every $r\ge 4$ the set $\Omega$ of points of $\cH_q$  is a complete $(k,q+1)$-arc with only four characters $0,1,2,q+1$ where $k=q^{2r}+1\pm q^{r+1}(q-1)$ according as $r$ is odd or even; see Theorem \ref{teoE}. For the rational BKS case,  the set $\Omega$ of the points in $\PG(2,q^r)$ is a $(k,q+1)$-arc with $k=q^r+1-\ha q(q-1)$ and characters  $0,1,2,\ha(q+1),\ha(q+3),q,q+1$. Furthermore, $\Omega$  is complete if and only if $r$ is even. If $r$ is odd, then the uncovered points by the $(q+1)$-secants to $\Omega$ are exactly the points in $\PG(2,q)$ not lying in $\Omega$. Adding those points to $\Omega$ produces a complete $(k,q+1)$-arc in $\PG(2,q^r)$, with $k=q^r+q$; see Theorem \ref{teoF}.

The above results do not hold true for $r=2$ and it remains open the case $r=3$ for the Hermitian curve, and the cases $r=3,4$ for the rational  BKS curve.

As a by product we also have the following two results of interest in the study of the Galois inverse problem.

Let $K=\bar{\mathbb{F}}_{q^{2r}}(m)$ and $L=K(u)$ where $u^{q+1}+u^qm^q+um-((ma-b)^q+ma-b)$ and $a^{q+1}+b^q+b\ne 0$. Then the geometric monodromy group of $L|K$ is
isomorphic to $\PGL(2,q)$, and the Galois closure $M$ of $L|K$ is
$M=\bar{\mathbb{F}}_{q^{2r}}(m,u,v,w)$ where
$$
\begin{cases}
u^{q+1}+u^qm^q+um-((ma-b)^q+ma-b)=0;\\
v^q+(u+m^q)v^{q-1}+u^q+m=0;\\
v+u+m^q-(u+m^q)w^{q-1}=0.
\end{cases}
$$

Let
$K=\bar{\mathbb{F}}_{q^{r}}(t)$ and $L=K(u)$ where $u^{q+1}+um^q+um-(b-2)(t-1)-\ha q+1$ and $b^{q+1}-(a^q+a)+(b^2-4a)^{(q+1)/2}\neq 0$. Then the geometric monodromy group of $L|K$ is
isomorphic to $\PGL(2,q)$, and the Galois closure $M$ of $L|K$ is
$M=\bar{\mathbb{F}}_{q^{r}}(m,u,v,w)$ where
$$
\begin{cases}
u^{q+1}+mu^q+mu-(b-2)(m-1)-\ha a+1=0, \\
v^q+(u+m)v^{q-1}+u^q+m=0, \\
v+u+m-(u+m)w^{q-1}=0.
\end{cases}
$$

\section{Outline of the proofs for the Hermitian case}
Some more notation is needed: $\cH_q$ denotes the (absolutely irreducible) Hermitian curve of homogeneous equation $Y^qZ+YZ^q+X^{q+1}=0$ defined over $\mathbb{F}_{q^2}$ and viewed as an (absolutely) irreducible curve in $\PG(2,q^{2r})$ for $r\ge 3$, and $\Omega$ stands for the set of all points of $\cH_q$ in $\PG(2,q^{2r})$ where $k=|\Omega|$ with $k=q^{2r}+1\pm q^{r+1}(q-1)$ according as $r$ is odd or even; see for instance \cite[Chapter 10]{HKT}. Then $\Omega$ is a $(k,q+1)$-arc in $\PG(2,q^{2r})$. For $r=1$, $\Omega$ is the classical unital, and hence it is a complete $(q^3+1,q+1)$-arc in $\PG(2,q^2)$. This does not hold true for $r=2$, as $\Omega$ in $\PG(2,q^4)$ is contained in $\PG(2,q^2)$ and hence no $(q+1)$-secant to $\Omega$ covers a point $P\in \PG(2,q^4)\setminus \PG(2,q^2)$ provided that $P$ is chosen on a tangent line to $\cH_q$.

To deal with the completeness problem in the general case, take any point $P\in \PG(2,q^{2r})$ not in $\PG(2,q^2)$. Since $PGU(3,q)$ leaves $\Omega$ invariant and preserves no line in $\PG(2,q^{2r})$, we may assume that $P$ is not a point at infinity. Therefore we use affine coordinates with $Z=0$ taken to be the line at infinity. Let $P=(a,b)$. If $\ell_m$ denotes the (non-vertical) line through $P$ with slope $m$, i.e. $Y=m(X-a)+b$, and
$$F(X)=X^{q+1}+X^q(a+m^q)+X(a^q+m)+b^q+b+a^{q+1} \in \mathbb{F}_{q^{2r}}[X]$$
then $\ell_m$ is a $(q+1)$-secant to $\Omega$ if and only if $F(X)$ has $q+1$ pairwise distinct roots in $\mathbb{F}_{q^{2r}}$. Now, take an algebraic closure $\mathbb{\bar{F}}$ of $\mathbb{F}_{q^2}$ containing $\mathbb{F}_{q^{2r}}$, and look at $F(X)$ as a polynomial with coefficients in the rational field $K=\mathbb{\bar{F}}(m)$. Two cases are distinguished according as $P$ lies in $\Omega$ or does not.

Assume first $P\notin \Omega$. Then $F(X)$ is an irreducible separable polynomial over $K$. Take a root $u$ of $F(X)$ in some overfield of $K$, and define $L=K(u)$ to be the algebraic extension of $K$ by adjoining $u$. The field extension $L|K$ is not Galois. The Galois closure $M$ of $L|K$ is the splitting field of $F(X)$ over $K$, and the associated Galois group $G=Gal(M|K)$ is the geometric monodromy group of $F(X)$ over $K$.
We prove in our case that $M=K(u,v,w)$. Our proof is based on Abhyankar's work \cite{abh}, especially on the concept of a twisted Abhyankar's derivative of a polynomial. More precisely, we show in Section \ref{abhy} that if the first Abhyankar derivative $f_1$ of $F$ is irreducible then the second Abhyankar derivative $f_2$ of $F$ splits into linear factors. By \cite[Section 2]{abh}, $G$ acts faithfully on the roots of $F$ as a sharply  $3$-transitive permutation group whose $2$-point stabilizer is cyclic. From Zassenhaus' classification of finite sharply $3$-transitive permutation groups \cite{Za}, $G=\PGL(2,q)$ follows. The missing piece in this argument, i.e. the irreducibility of $f_1$, is proven in Section \ref{vWt} where we rely  on a classical theorem of van der Waerden \cite{vW}.
The result $G= \PGL(2,q)$ is quite a surprising  since in most cases $2$-transitive geometric monodromy groups are either the symmetric group or the alternating group.

In our case $G\cong \PGL(2,q)$. To prove it we exploit Abhyankar's work \cite{abh} and a theorem of van der Waerden \cite{vW} together to show that $G$ acts faithfully on the roots of $F(X)$ as a $2$-transitive permutation group. A closer look at the action of $G$ allows to prove that the $2$-point stabilizer is cyclic, and then $G\cong \PGL(2,q)$ follows from a characterization of $\PGL(2,q)$ as the unique sharply $3$-transitive group with cyclic $2$-point stabilizer.
This is quite a surprising result since in most cases $2$-transitive geometric monodromy groups are either the symmetric group or the alternating group.

The next step is to show that the ramified places in the Galois extension $M|K$  are as many as $(q+1)^2$. From this we deduce that $G$ has $q+1$ short orbits on the set of places of $M$ and that it acts on each short orbit as $\PGL(2,q)$ in its $3$-transitive permutation representation.

It turns out that the point $P$ is covered by at least one (non-vertical) line $\ell_m$ if and only if $M$ has at least one $\mathbb{F}_{q^{2r}}$-rational place unramified in the Galois extension $M|K$. Using Serre's ramification theory \cite{serre1979}, see also \cite[Section 11.9]{HKT}, we are able to compute the genus $\mathfrak{g}(M)$ of $M$. Actually, $\mathfrak{g}(M)$ only depends on $q$, as
$2\mathfrak{g}(M)-2=q^4-q^2-2q-2.$ Now, if
$$q^{2r}+1> 2\mathfrak{g} q^r+(q+1)^2>q^{r+4}-q^{r+2}-2q^{r+1}+q^2+2q+1,$$
then the Hasse-Weil lower bound ensures the existence of $m\in \mathbb{F}_{q^{2r}}$ such that the polynomial
$F(X)$ has $q+1$ pairwise distinct roots over $\mathbb{F}_{q^2}$. Therefore, $P$ is covered by a $(q+1)$-secant to $\Omega$.  Thus $r=3$ remains the only open case. A Magma aided search shows that for $q=r=3$ $\Omega$ is complete.

The case $P\in \Omega$ is treated analogously. Let $P=P(a,b)$ with $b^q+b+a^{q+1}=0$, and
$$F(X)=X^q+X^{q-1}(a+t^q)+(a^q+t) \in \mathbb{F}_{q^{2r}}[X].$$
Then $G$ acts faithfully on the roots of $F(X)$ as a sharply $2$-transitive permutation group, and $G\cong \AGL(1,q)$. Furthermore, $G$ fixes a place of $M$ and has a unique non-trivial short orbit of size $q$. From this,
$\mathfrak{g}(M)=\ha q(q-1)^2$ follows. If
$$q^{2r}+1> 2\mathfrak{g} q^r+q+1>q^{r+3}-2q^{r+2}+q^{r+1}+q+1,$$
then the Hasse-Weil lower bound yields that $P$ lies on a $(q+1)$-secant to $\Omega$. This is indeed the case since as $r\ge 3$ has been assumed.

\section{Preliminaries on absolutely irreducibility of polynomials and plane curves}
Let
\begin{equation}
\label{eqB8dic22}
P(X)=X^{q+1}+eX^q+aX+b\in \mathbb{F}_{q^s}[X].
\end{equation}
As Bluher \cite{bluh} pointed out if $ea \neq b$ and $a\neq e^q$ then the substitution of $X$ by $(e^{q+1}-b)/(a-e^q)X-e$ brings $P(X)$ to the form $X^{q+1}-BX+B$ where $B=(a-e^q)^{q+1}/(b-ea)^q\in  \mathbb{F}_{q^s}^*$. In particular, it has no multiple roots in $\mathbb{F}$. Also, she proved that if $X^{q+1}-BX+B$ has at least three (distinct roots) then all the $q+1$ roots are in  $\mathbb{F}_{q^s}$; see \cite[Theorem 4.3]{bluh}.

\begin{result}
\label{bluherA} Set either
$$P(X)=
\begin{cases}
X^{q+1}+m^qX^q+mX-((ma-b)^q+ma-b)\in \mathbb{F}_{q^s}[X],\\
m^{q+1}+(ma-b)^q+(ma-b)\neq 0,\\
\end{cases}
$$
or
$$
P(X)=
\begin{cases}
2mX^{q+1}+(2m-1)X^q+(2m-1)X+m(2-a)+b-2\in \mathbb{F}_{q^s}[X],\\
 2am^2-2mb+1\neq 0.
\end{cases}
$$
Then $P(X)$ has no multiple roots. Furthermore, if $P(X)$ has at least three roots in $\mathbb{F}_{q^s}$ then all its $q+1$ roots are in  $\mathbb{F}_{q^s}$.
\end{result}

\begin{lemma}
\label{lem5dic} Let $a,b\in \mathbb{F}$ such that $\mathbb{F}$ has no element $t$\, for which $a=2(t+1)^{q+1}$ and $b=2+t^q+t$.
Then the plane curve $\cC$ with affine equation
$$F(U,V)=V(2(U+1)^{q+1}-a)-U^q-U-2+b=0$$
is irreducible.
\end{lemma}
\begin{proof} The point at infinity $V_\infty$ is an ordinary singular point of $\cC$ with multiplicity $q+1$. The tangent lines $\ell_i$ to $\cC$ at $V_\infty$ have equations $U-u_i=0$ where $(u_i+1)^{q+1}=\ha a$ and $i=1,2,\ldots, q+1$. None of them is a linear component of $\cC$ as $b\ne 2+u_i^q+u_i$. Moreover $I(V_\infty, \cC\cap \ell_i)=q+2$. If $\cC$ is reducible, Segre's criterium, \cite[Lemma 8]{BS} applies to each  $\ell_i$. Therefore, $P(U)=\prod_{i=1}^{q+1} (U-u_i)$ divides $F(U,V)$. Since $P(U)=\prod_{i=1}^{q+1} (U+1-(u_i+1))=(U+1)^{q+1}-\ha a$, this yields $F(U,V)=(2(U+1)^{q+1}-a)F_1(U,V)$ with $\deg(F_1(U,V))=1$. Thus  $(2(U+1)^{q+1}-a)(V-F_1(U,V))-(U^q+U+(b-2))$ would be the zero polynomial, a contradiction. Therefore, $\cC$ is irreducible.
\end{proof}
\begin{lemma}
\label{lem7dic} Let\, $t\in \mathbb{F}\setminus \mathbb{F}_{q^2}$. Then the plane curve $\cC$ with affine equation
$$F(U,V)=V(U^{q-1}+1)+U^q+t U^{q-1}+t^q=0$$
is irreducible.
\end{lemma}
\begin{proof} For $i=1,\ldots, q-1$, let $\ell_i$ denote the line of equation $U-u_i=0$ with $u_i^{q-1}+1=0$. If $\ell_i$ is a component of $\cC$  then $u_i^q+tu_i^{q-1}+t^q=0$ and hence $t^q-t=u_i$. But then $(t^q-t)^{q-1}=-1$ whence $(t^q-t)^q=-(t^q-t)$ from which $t^{q^2}=t$, that is, $t\in \mathbb{F}_{q^2}$. Therefore, $\ell_i$ is not a component of $\cC$. Therefore, the point at infinity $V_\infty$ is an ordinary singular point of $\cC$ with multiplicity $q-1$, and the tangent lines $\ell_i$ to $\cC$ at $V_\infty$ have equations $U-u_i=0$ where $u_i^{q-1}=-1$ for $i=1,2,\ldots q-1$. Now, we argue as in the proof of Lemma \ref{lem5dic}. From Segre's criterium, $P(U)=\prod_{i=1}^{q-1}(U-u_i)=U^{q-1}+1$ divides $F(U,V)$. Thus $F(U,V)=(U^{q-1}+1)F_1(U,V)$  with $\deg(F_1(U,V))=1$ whence
the claim follows.
\end{proof}
\label{bfg}
\section{$(k,n)$-arcs arising from the Hermitian curve and the rational BKS curve}
\subsection{The Hermitian curve and its geometry}
\label{hhc}
Let $\cH_q$ denote the Hermitian curve given in its canonical form of affine equation
\begin{equation}
\label{eq3nov24}
X^{q+1}+Y^q+Y=0.
\end{equation}
The properties of $\cH_q$ pertinent to the present paper are: (i) $\cH_q$ is non-singular of genus $\mathfrak{g}=\ha q(q-1)$, (ii) any line $\ell$ of  $\PG(2,q^2)$ either meets $\cH_q$ in $q+1$ pairwise distinct points all lying in $\PG(2,q^2)$, or is a tangent to $\cH_q$ at a point $P\in \PG(2,q^2)$ and $I(P,\cH_q\cap \ell)=q+1$, in particular, each common point of $\ell$ with $\cH_q$ lies in $\PG(2,q^2)$, (iii) $\Omega$ is a complete $(q^3+1,q+1)$-arc with two characters namely $1$ and $q+1$, (iv) the subgroup of $\PGL(3,q^2)$ which leaves $\cH_q$ invariant is the projective unitary group $PGU(3,q)$, (v) $\cH_q$ is an $\mathbb{F}_{q^2}$-maximal curve and the set of its points lying in $\PG(2,q^2)$ has size $q^3+1$.

From now on, we focus on the case $r\ge 3$. The set $\Omega$ consisting of all points of $\cH_q$ lying in $\PG(2,q^{2r}$ is equal to $q^{2r}+1\pm q^{n+1}(q-1)$ according as $r$ is odd or even \cite[Chapter 10]{HKT}. Fix a point $Q\in \PG(2,q^{2r})$ not in $\PG(2,q^2)$. Since $\PGL(3,q)$ does not leave the infinite line $Z=0$ invariant, we may assume $Q\notin \ell_\infty$ and use affine coordinates. Then $Q=Q(a,b)$ with $a,b\in \mathbb{F}_{q^{2r}}$. For any line $\ell$ through $Q$ we determine its common points with $\cH_q$. If the vertical line coincides with $\ell$ then $\cH_q\cap \ell$ comprises $Y_\infty$ together with the points $P(a,t)$ such that $t$ is the root of the polynomial $F(Y)=Y^q+Y+a^{q+1}=0$. If one of the roots belongs to $\mathbb{F}_{q^{2r}}$ then all do. Also, $F(Y)$ is separable and hence  $\ell$ is not a tangent to $\cH_q$. Now, let $\ell$ be a line through $Q$ of equation $Y=m(X-a)+b$. Then the common points of $\ell$ and $\cH_q$ are the points $P(\xi,\eta)$ such that $\xi$ is a root of the polynomial
\begin{equation}
\label{eq241122}
F(X)=X^{q+1}+m^qX^q+mX-((ma-b)^q+ma-b).
\end{equation}
which can also be written as
\begin{equation}
\label{eq2411A22}
F(X)=(X^q+m)(X+m^q)-(m^{q+1}+(ma-b)^q+(ma-b)).
\end{equation}
Then $\xi$ is a multiple root if and only if $\xi$ is a root of the polynomial $dF/dX=X^q+m$ as well, that is, $\xi^q+m=0$. From (\ref{eq2411A22}) this occurs if and only if $m^{q+1}+(ma-b)^q+ma-b=0$. The polynomial $G(T)= T^{q+1}+(Ta-b)^q+Ta-b$ has a multiple root if and only if $b^q+b+a^{q+1}=0$, that is, $Q(a,b)\in \cH_q$. Therefore two cases arise. If $Q(a,b)\notin \cH_q$ then there exist $q+1$ lines $\ell_m$ which are tangent to $\cH_q$ and the tangency point of $\ell_m$ is $P_m=(-\sqrt[q]{m},m(\xi-a)+b)$ with $\xi^q=-m$.  Furthermore, $\ell_m$ also meets $\cH_q$ in the point $R_m=(-m^q,m(\xi-a)+b)$ with $\xi=-m^q$. It is possible that $P_m=R_m$ and this occurs when $\xi^{q^2}-\xi=0$, that is, $\xi\in \mathbb{F}_{q^2}$. In this case $I(P_m,\cH_q\cap \ell_m)=q+1$. Otherwise, $P_m\neq R_m$ where $I(P_m,\cH_q\cap \ell_m)=q$ and $I(R_m,\cH_q\cap \ell_m)=1$. If $Q\in \cH_q$ then $\cH_q$ has a unique tangent $\ell$ at $Q$ where $\ell$ has equation $Y=-a^q(X-a)+b=-a^qX-b^q$. Furthermore, $R=(a^{q^2},b^{q^2})$ is another common point of $\cH_q$. Since $Q\neq R$ by our assumption, $I(Q,\cH_q\cap \ell)=q$ and $I(R,\cH_q\cap \ell)$=1. Therefore, any tangent to $\cH_q$ through $Q$ meets $\Omega$ in either one or two points, and in the former case the tangency point is in $\PG(2,q^2)$.   Take a line $\ell$ be through $Q$ other than tangents $\ell_i$. Then $\ell\cap \cH_q$ consists of $q+1$ pairwise distinct points lying in some extension of $\PG(2,q^{2r})$. Furthermore, if three of them are in $\PG(2,q^{2r})$ then each of them is in $\PG(2,q^{2r})$. In fact, if three roots of the polynomial $F(T)$ given in (\ref{eq241122}) belong to $\mathbb{F}_{q^{2r}}$ then all do; see Result \ref{bluherA}.

The above results show that the set $\Omega$ of all points of $\cH_q$ in $\PG(2,q^{2r})$ with $r\ge 2$ is a $(k,q+1)$-arc with characters $0,1,2,q+1$ where $k=q^{2r}+1\pm q^{r+1}(q-1)$ according as $r$ is odd or even.

\subsection{The rational BKS curve and its geometry}
\label{BKSs}
For $q$ odd, let $\cC$ be the curve of affine equation
\begin{equation}
\label{eq3jul22}
Y^{q+1}-(X^q+X)+(Y^2-2X)^{(q+1)/2}=0.
\end{equation}
For the known properties of $\cC$; see \cite[Proposition 4.24]{BKS}: $\cC$  has as many as $\ha q(q-1)$ singular points each of them being a node (ordinary doubly point) lying in $\PG(2,q)$ with both tangents defined over $\mathbb{F}_{q^2}$. The singular points of $\cC$ are exactly the internal points to the conic $\cC^2$ of affine equation $Y^2-2X=0$. The intersection of $\cC$ with a tangent line $\ell$ at a singular point $P$ of $\cC$ collapses into $P$. More precisely, if $\gamma$ is a branch of $\cC$ centered in $P$ then the intersection multiplicity $I(P,\gamma\cap \ell)=q$
whereas $I(P,\delta\cap \ell)=1$ for the other branch $\delta$ of $\cC$ centered in $P$. The $q+1$ points of $\cC^2$ in $\PG(2,q)$ are also points of $\cC$. The intersection of $\cC$ with the tangent line $\ell$ at a point $P\in\cC^2$ of $\cC$ also collapses into $P$, that is, $I(P,\cC\cap \ell)=q+1$. The singular points of $\cC$ together with the points of $\cC^2$ lying in $\PG(2,q)$ form the set of size $\ha q(q-1)+q+1$
consisting of all points of $\cC$ lying in $\PG(2,q^2)$. The projective closure $\cD$ of $\cC$ is invariant under the action of a subgroup $G\cong \PGL(2,q)$ of $\PGL(3,q)$ which acts on $\cC\cap \cC^2$ as $\PGL(2,q)$ in its unique $3$-transitive permutation representation.
Since $\cC$ has degree $q+1$ and possesses $\ha q(q-1)$ nodes, the genus of $\cC$ equals zero and hence $\cC$ is a rational curve. Thus $\cC$ can be parametrized by a variable $t$ over $\bar{\mathbb{F}}_q$. More precisely,
$\cC$  consists of the points
\begin{equation}
\label{eqpoints}
P(t)=(2(t+1)^{q+1},2+t+t^q),\,\, t\in \bar{\mathbb{F}}_q^*\cup \{\infty\}
\end{equation}
where $\infty$ stands for the parameter of the point at infinity $P_\infty=(1:0:0)$. In this parametrization, $P(t)\in \PG(2,q)$ if and only if $t\in \mathbb{F}_{q^2}\cup\{\infty\}$ where either $t\in \mathbb{F}_q\cup \{\infty\}$ or  $t\in \mathbb{F}_{q^2}\setminus \mathbb{F}_q$ holds according as $P(t)\in \cC^2\cap \PG(2,q)$, or $P(t)$ is an internal point to $\cC^2$ in $\PG(2,q)$.
%Furthermore, $P(t)\in \PG(2,q)$ if and only if either $t\in \mathbb{F}_q$ and $P(t)\in \cC^2$, or $t\in \mathbb{F}_{q^2}\setminus \mathbb{F}_q$ and $P(t)\in \PG(2,q)$ is an internal point to $\cC^2$.
In the latter case,
case $P(t)=P(t^q)$. For the other points, $P(t)=P(t')$ only occurs for $t=t'$, and $P(t)\in \PG(2,q^r)$ with $r\ge 3$ if and only $t\in \mathbb{F}_{q^r}$.

If a line $\ell$ is defined over $\mathbb{F}_q$ then the above properties of $\cC$ completely determine $\ell \cap \cC$: (i) If $\ell$ is an external line to $\cC^2$ in $\PG(2,q)$ then $\ell\cap \cC$ consists of $\ha (q+1)$ points each being a double points of $\cC$; (ii) if $\ell$ is a chord of $\cC^2$ in $\PG(2,q)$ with $\ell\cap \cC^2=\{R_1,R_2\}$ then $\ell\cap \cC$ comprises $R_1,R_2$ and $\ha (q-1)$ points each being a double points of $\cC$;(iii) if $\ell$ is a tangent line to $\cC^2$ in $\PG(2,q)$ with $\ell\cap \cC^2=\{R\}$ then the intersection collapses into the point $R$. This shows that the set $\Omega$ consisting of all points of $\cC$ lying in $\PG(2,q)$ is a
$(k,n)$-arc with $k=\ha q(q-1)+q+1=\ha(q^2+q+2), n=\ha (q+3)$ and three characters $1,\ha(q+1),\ha(q+3)$. In particular, $\Omega$ is not a $(k,q+1)$-arc although it arises from a curve of degree $q+1$. Moreover, since no point of $\PG(2,q^2)\setminus \PG(2,q)$ belongs to $\cC$, the set consisting of all points of $\cC$ lying in $\PG(2,q)$ will have the same three characters $1,\ha(q+1),\ha(q+3)$.

In $\PG(2,q^r)$ with $r\ge 1$, fix a point $Q$ other than those of $\cC$ lying in $\PG(2,q)$.
We are going to study the possible intersection of $\cC$ with  a line $\ell$ through $Q$.

If $\ell$ has equation $X=a$ then $\ell$ cuts out on $\cC$ the points of parameters $t$ that are the solutions of the equation $(t+1)^{q+1}+\ha a=0$. If $\ell$ has equation $Y=b$ then the parameters of the points cut out by $\ell$ on $\cC$ are $\infty$ together with the solutions of the equation $t^q+t+2=b$. If $\ell$ has equation $Y=\ha(X-a)+b$ then the points cut out by $\ell$ on $\cC$ have parameters $t$ satisfying the equation $t^{q+1}=1+\ha a-b$. In all three cases, the solutions of the equations are pairwise distinct and define distinct points of $\cC$ unless one of the following cases occurs: either, each of these points is counted twice and $\ell$ is a line of $\PG(2,q)$, or $a=0$, or $b=\ha a +1$ and $t=0$ is the unique solution. In the latter two cases, $\ell$ is a line of $\PG(2,q)$, as well, since $\ell$ has equation $X=0$ and $Y=\ha X+1$, respectively.

From now on, the line $\ell$ through $Q(a,b)$ is assumed to be distinct from the lines of $\PG(2,q)$.

To determine when $\ell$ is a tangent to $\cC$, we may assume that $\ell$ has slop $m$ with $m\in \mathbb{F}$, and $m\ne 0,\ha$. %Moreover, we replace $t$ by $t^{-1}$ in the parametrization (\ref{eqpoints}) of $\cC$.
Then the parameters $t$ of the points of $\cC$ cut out by $\ell_m$ are the roots of the polynomial $m(2(T+1)^{q+1}-a)-2-T^q-T+b$ that is of
\begin{equation}
\label{eq3jul2022C} F(T)=2mT^{q+1}+(2m-1)T^q+(2m-1)T+m(2-a)+b-2,
\end{equation}
Then $t$ is a multiple root of $f(T)$ if and only if $t$ is also a root of the polynomial $dF(T)/dT=2m(T+1)^q-1$. This together with $F(t)=0$ yield $(2m-1)t^q+m(2-a)+b-2=0$. Eliminating $t^q$ from the equations $2m(t+1)^q=1$ and $(2m-1)t^q=-b-2m-ma+2$ gives
$2am^2-2mb+1=0$. Also, if $m$ satisfies  $2am^2-2mb+1=0$ then $$F(T)=\frac{1}{2m}(2mT+(2m-1))(2mT^q+(2m-1)),$$ and the converse also holds. From this, if $b^2-2a \neq 0$ then through $Q=Q(a,b)$ there are exactly two tangents to $\cC$: For $a\ne 0$, the lines $\ell_i$ of equation $Y=m_i(X-a)+b$ with
%$$m_{1,2}=\frac{b\pm \sqrt{b^2-2a}}{2a}.$$
$m_i=(b\pm \sqrt{b^2-2a})/(2a)$ while for $a=0$ the line $m_1$  of equation $Y=\ha b^{-1}X+b$ and the line $r_2$ of equation $X=0$. For $b^2-2a=0$, that is, for $Q\in\cC^2\setminus \cC$, we have $m_1=m_2$ and there is a unique tangent $\ell_1$ to $\cC$ where $\ell_1$ has equation $Y=b/(2a)(X-a)+b$ and hence is the tangent to $\cC^2$ at $Q$.
Also, for $i=1,2$, $\ell_i$ has exactly two common points with $\cC$, namely the points $P_i,Q_i$ whose parameters are $(1-2m_i)/(2m_i)$ and $\sqrt[q]{(1-2m_i)/(2m_i)}$. Here, $I(P_i,\cC\cap \ell_i)=1$ and  $I(Q_i,\cC\cap \ell_i)=q$.
unless $P_i=Q_i$ and $I(Q_i,\cC\cap \ell_i)=q+1$. In the latter exceptional case, $m_i\in \mathbb{F}_q$. Furthermore, $\ell_i$ is the tangent to $\cC$ at $Q_i$.
This result holds true for $r_2$ which is the tangent to $\cC$ at $O=(0,0)$ where $I(O,\cC \cap r_2)=q+1$. From now on set $r_2=m_2$. With this notation, $m_1,m_2$ are the tangents to $\cC$ through $Q(a,b)$ where $m_1=m_2$ if and only if $Q\in \cC^2\setminus \cC$. In any case, $m_i\in \mathbb{F}_{q^{2r}}$ where $m_i\in \mathbb{F}_{q^r}$ occurs whenever $b^2-2a$ is a square in $\mathbb{F}_{q^r}$. Accordingly, both tangency points $P_i,Q_i$ are in $\PG(2,q^{2r})$ or in $\PG(2,q^r)$. Here, $P_i=Q_i$ is only possible when $m_i\in \mathbb{F}_q$. In particular, both $P_1=Q_1$ and $P_2=Q_2$ hold if and only if $Q$ is an external point to $\cC^2$ in $\PG(2,q)$ and the tangents through $Q$ to $\cC$ and those to $\cC^2$ coincide. Both tangents are lines of $\PG(2,q)$, and this case has already been discussed before.

Finally, let $\ell$ be a line through $Q$ other than the tangents $\ell_i$ and $r_i$. If $\ell$ passes through a (unique) singular point $R\in \PG(2,q)$  of $\cC$ then $\ell\cap \cC$ consists of $R$ together with $q-1$ pairwise distinct points lying in $\PG(2,q^r)$ or in some extension of $\PG(2,q^r)$. Otherwise, $\ell\cap \cC$ consists of $q+1$ pairwise distinct points lying in $\PG(2,q^r)$ or some extension of $\PG(2,q^r)$. As in the Hermitian case, if three of them are in $\PG(2,q^r)$ then each of them is in $\PG(2,q^r)$; see Result \ref{bluherA}.  Also, since each internal point to $\cC$ corresponds to two parameters, $\Omega$ has size $q^r+1-\ha q(q-1)$

From the above results, the set $\Omega$ of all points of $\cC$ in $\PG(2,q^r)$ with $r\ge 2$ is a $(q^r+1-\ha q(q-1),q+1)$-arc with characters $0,1,2,\ha(q+1),\ha(q+3),q,q+1$.

%\section{Background from Galois theory} We use standard notation and terminology from Function field theory:
%\label{gth}
\section{Abhyankar's work}
\label{abhy}
An important tool for the study of the action of the Galois group on the roots of its defining polynomial, given by explicit equation, is Abhyankar's skew derivative introduced in  \cite{abh}. Let
\begin{equation}
\label{eq241122B}
f=f(T)=T^m+a_1T^{m-1}+\ldots+ a_m
\end{equation}
 be a separable, monic polynomial in the indeterminate $T$ with coefficients $a_i$ in a field $K$. The splitting field $M$ of $f(T)$ is generated by $K$ together with the roots $\alpha_1,\alpha_2,\ldots, \alpha_m$ of $f(T)$, and the Galois group $Gal(M|K)$ consists of all $K$-automorphisms of $M$. Furthermore, $Gal(M|K)$ acts faithfully on the set $\Delta=\{ \alpha_1,\alpha_2,\ldots, \alpha_m\}$, and it can be viewed as a permutation group on $\Delta$ named the Galois group $Gal(f,K)$ of $f$ (over $K$). The group $Gal(f,K)$ is transitive on $\Delta$ if and only if $f(T)$ is irreducible (over $K$).

 From now on assume that $F(T)$ is irreducible.  Following \cite[Section 4]{abh}, we  ``throw away'' a root of $f(T)$, say $\alpha_1$, and get
 \begin{equation}
\label{eq241122C}
f_1=f_1(T)=\frac{f(T)}{T-\alpha_1}=T^{m-1}+b_1T^{m-1}+\ldots+ b_m \in K(\alpha_1)[T].
\end{equation}
Then  $f_1(T)$ is irreducible over $K(\alpha_1)$ if and only if $Gal(f,K)$ is $2$-transitive on $\Delta$. As Abhyankar stressed in \cite[Section 4]{abh}, it does not matter which root of $f(T)$ we throw away; for instance, the irreducibility of $f_1(T)$ over $K(\alpha_1)$
and, up to isomorphism, the Galois group $Gal(f_1, K(\alpha_1))$ are independent of which root we call $\alpha_1$. Likewise, by throwing away $s$ roots of $f$ we get
 \begin{equation}
\label{eq241122C1}
f_s=f_s(T)=\frac{f(T)}{(T-\alpha_1)\cdots (T-\alpha_s)}=T^{m-s}+d_1T^{m-s-1}+\ldots+ d_{m-s} \in K(\alpha_1,\ldots,\alpha_s)[T],
\end{equation}
and $f_s$ is irreducible over $K(\alpha_1,\ldots,\alpha_s)$ for $1\le s \le m$ if and only if $Gal(f,K)$ is $s$-transitive on $\Delta$.

 For any polynomial $\Theta = \Theta(T)$ in an indeterminate $T$ with coefficients
in a field $L$ and for any element $\beta \in L$,
$$\frac{\Theta(T + \beta) - \Theta(\beta)}{T}$$
is the first Abhyankar's twisted $T$-derivative of $T$  at $\beta$; see \cite[Art.33]{abh}.  With this definition, the following result is stated in \cite[pg.\,93]{abh}: If $f= f(T)$ is a nonconstant monic
irreducible polynomial in an indeterminate $T$ with coefficients in a field $K$
such that$f$ has no multiple root in any overfield of $K$, and if $\alpha$ is a root
of $f$ in some overfield of $K$, then by letting $f' = f'(T)$ to be the twisted
derivative of $f$ at $\alpha$ we have that the Galois group $Gal(f', K(\alpha))$ is the
one-point stabilizer of the Galois group $Gal(f,K)$.
\subsection{The Hermitian case} Let $K=\mathbb{F}(m)$ be the rational field over an algebraically closed field $\mathbb{F}$ of positive characteristic $p$. Fix a power $q$ of $p$.
Consider the polynomial
\begin{equation}
\label{eq1020322}
f=f(T)=T^{q+1}+m^qT^q+mT-((ma-b)^q+ma-b)\in \mathbb{F}(m)[T].
\end{equation}
\subsubsection{Case  $a,b\in \mathbb{F}$ with $a^{q+1}+b^q+b\neq 0$}
\label{subs1}
\begin{lemma}
\label{lemA271122}
$f(T)$ is irreducible over $K$.
\end{lemma}
\begin{proof} In the Hermitian function field $H_q=\mathbb{F}(x,y)$ with $y^q+y-x^{q+1}=0$, let $\varphi$ be rational map defined by $\varphi(x)=x,\varphi(y)=(y-b)/(x-a)$. Clearly, $\varphi$ is birational as $y=\varphi(y)(x-a)+b$. Set $m=\varphi(y)$. Let $g=g(X,Y)$ be a minimal polynomial of $x$ and $m$. Then $H_q$ coincides with the function field $U=U(x,m)$ where $g(x,m)=0$. Now, regard (\ref{eq1020322}) as a polynomial $f(X,Y)$ over $K$ with $f(X,Y)=X^{q+1}+Y^qX^q+XY-((Ya-b)^q+Ya-b)$. Choose an irreducible non-constant factor $h(X,Y)\in \mathbb{F}[X,Y]$ of $f(X,Y)$. Take $\xi,\mu \in \mathbb{F}$ such that $h(\xi,\mu)=0$. Then $f(\xi,\mu)=0$ whence $\eta^q+\eta-\xi^{q+1}=0$ for $\eta=\mu(\xi-a)+b$. This yields $g(\xi,\mu)=0$. From Study's theorem \cite[Theorem 2.10]{HKT}, $h$ divides $g$. Since $g$ is irreducible, this yields $h=g$. Assume on the contrary that (\ref{eq1020322}) reducible.  Then $f(X,Y)=c g(X,Y)^i$ for $i\ge 2$ and a non-zero constant $c$. Since $\deg(f(X,Y))=2q$ this yields $\deg(g(X,Y))\le q$. On the other hand, $H_q=U$ implies that these function fields have the same genus $\ha q(q-1)$. But then $\deg(U)\ge q+1$, a contradiction.
\end{proof}
We point out that the curve $\cC$ of affine equation $f(X,Y)=0$ introduced in the proof has two singular points, namely the points at infinity of the $X$ and $Y$ axes. Actually, each other point of $\cC$ is non-singular. In fact, $f_X(\xi,\eta)=0$ and $f_Y(\xi,\eta)=0$ yield $\eta=-\xi^q$ and $\xi=a$, but $f(\xi,\eta)=U(a,-a^q)=-(a^{q+1}+b^q+b)\neq 0$.

Let $u$ be a root of $f(T)$ in some extension of $K$ and put $L=K(u)$. A straightforward computation shows that the first Abhyankar's skew derivative of $f$ at $u$ is
\begin{equation}
\label{eqA251122}
f_1=f_1(T)=\frac{f(T+u)-f(u)}{T}=T^q+(u+m^q)T^{q-1}+(m+u^q)\in L[X].
\end{equation}

Now we compute the second Abhyankar's twisted derivative of $f(T)$, i.e. the Abhyankar's twisted derivative of $f_1(T)$ at any $v$ which is a root of $f(T)$ different from $u$.
\begin{equation}
\label{eq3251122}
\frac{f_1(T)-f_1(v)}{T}=%T^{q-1}+(u+m^q)\frac{(T+v)^{q-1}-v^{q-1}}{T}=
T^{q-1}+(u+m^q)v^{q-1}\frac{(T/v+1)^{q-1}-1}{T}.
\end{equation}
Let $\Psi(U)$ be the polynomial whose roots are those of $f(T)$ different from $u$. This means replacing $T$ with $vU$ but preserving the splitting field $M$. Therefore,

%Replacing $T$ with $vU$ and dividing by $v^{q-2}$ give
$$ \Psi(U)=v^{q-1}U^{q-1}+(u+m^q)v^{q-2}\frac{(U+1)^{q-1}-1}{U}=$$
whence
$\Psi(U)=v^{q-2}(vU^{q-1}+(u+m^q)((U+1)^{q-2}+(U+1)^{q-3}+\ldots+1).$
We omit the factor $v^{q-2}$. Then
$$
\begin{array}{lll}
\Psi(U)=vU^{q-1}+(u+m^q)((U+1)^{q-2}+(U+1)^{q-3}+\ldots+1)= \\
vU^{q-1}-(u+m^q)(U+1)^{q-1}+(u+m^q)((U+1)^{q-1}+(U+1)^{q-2}+\ldots+1)= \\
vU^{q-1}-(u+m^q)(U+1)^{q-1}+(u+m^q)((U+1)^q-1)/U)=\\
vU^{q-1}-(u+m^q)(U+1)^{q-1}+(u+m^q)U^{q-1}=\\
(v+u+m^q)U^{q-1}-(u+m^q)(U+1)^{q-1}.
\end{array}
$$
Let $\Gamma(V)=(v+u+m^q)-(u+m^q)(V+1)^{q-1}$ be the polynomial obtained by reciprocating the roots of $\Delta(U)$, i.e. replacing $U$ with $V=U^{-1}$. This does not alter the splitting field $M$. Finally, let $\Phi(W)$ be the polynomial obtained by adding $1$ to the roots of $\Gamma(V)$, i.e. replacing
$V$ by $W=V+1$. Again, $M$ is left invariant. Then
\begin{equation}
\label{eq251122}
\Phi(W)=(v+u+m^q)-(u+m^q)W^{q-1}.
\end{equation}
Therefore the second Abhyankar's skew derivative of $f(T)$ at $v$ is
$$f_2(T)=v^{q-2}((v+u+m^q)T^{q-1}-(u+m^q)(T+v)^{q-1}).$$
We show that if $f_1(T)$ is irreducible over $K(u)$ then $f_2(T)$ is irreducible over $K(u,v)$. Assume on the contrary that $\Phi(W)$ is reducible. Then there exists $\rho\in K(u,v)$ such that $\rho^{q-1}(u+m^q)=v+u+m^q$ otherwise the $\Phi(W)$ defines a cyclic Kummer extension of $K(u,v)$ of degree $q-1$; see for instance \cite[Appendix A.13]{sti}. Then $v^q+(u+m^q)v^{q-1}+m+u^q=(v\rho)^{q-1}(u+m^q)+m+u^q=0$ where $(u+m^q)(m+u^q)\neq 0$.
Since $v\rho\in K(u,v)$, and $K(u,v):K(u)]=q$ this yields that $v\rho \in K(u)$. From Kummer's theory, there exists $\tau\in K(u)$ such that $-(m+u^q)/(m^q+u)=\tau^{q-1}$. Choose a root $m_0\in \mathbb{F}$ of the polynomial
$Y^{q+1}+(Ya-b)^q+Ya-b$ and then $u_0\in \mathbb{F}$ such that $u_0^q=-m_0$. Then $P(m_0,u_0)$ is a point of the curve $\cU$ with function field is $\mathbb{F}(m,u)$. As we have already shown that $\cU$ is non-singular. Therefore, $P(m_0,u_0)$ is a non singular point and the tangent to $\cU$ at $P(m_0,u_0)$ is the line of equation $U=m_0$. In particular, $m-m_0$ is a local parameter at $P(m_0,u_0)$. Therefore the valuation $v_{P(m_0,u_0)}(u^q+m)=1$. Furthermore, $u_0+m_0^q\neq 0$ otherwise $m\in \mathbb{F}_{q^2}$. Thus $v_{P(m_0,u_0)}(u+m^q)=0$. Therefore
$$v_{P(m_0,u_0)}\Big(\frac{u^q+m}{u+m^q}\Big)=1$$
which contradicts $-(m+u^q)/(m^q+u)=\tau^{q-1}$ with $\tau\in K(u)=\mathbb{F}(m,u)$.

Let $w$ be a root of $F(t)$ other than $u$ and $v$.
The third  Abhyankar's skew derivative $F_3(T)$ of $f(T)$ at $w$ is obtained by means of the first Abhyankar's skew derivative $\Phi_1=\Phi_1(W)$ of $\Phi(W)$ at $w$.
From (\ref{eq251122}),
$$\Phi_1(W)=\frac{\Phi(W+w)-\Phi(w)}{W}=-(u+m^q)w^{q-2}\frac{(W/w+1)^{q-1}-1}{W/w}$$
whence
$$\Phi_1(W)=-\frac{u+m^q}{w}\prod_{i=1}^{q-1}(W-\theta^i w)$$
for a primitive ($q-1$)-root $\theta$ of unity in $\mathbb{F}$.
This shows that $\Phi_1(W)$ (and hence $f_2(T)$) is a completely reducible polynomial over $K(u,v,w)=\mathbb{F}(m,u,v,w)$.
Moreover, $Gal(M|K)$ is a sharply $3$-transitive group on $\Delta$ such that the $2$-point stabilizer is cyclic.
From Zassenhaus' theorem \cite{Za}, $Gal(M/K)\cong \PGL(2,q)$, and $Gal(M/K)$ acts on $\Delta$ as $\PGL(2,q)$ on the projective line over $\mathbb{F}_q$.
\begin{theorem}
\label{teoA} Suppose that $f_1(T)$ is irreducible over $K=\mathbb{F}(m)$. Then $f_2(T)$ is irreducible over $K(u)$, and $M$ coincides with $K(u,v,w)=\mathbb{F}(m,u,v,w)$ where the extension $K(u,v,w)|K(u,v)$ is cyclic of degree $q-1$, and $M=\mathbb{F}(m,u,v,w)$ with
with
\begin{equation}
\label{eq261122}
\begin{cases}
u^{q+1}+m^qu^q+mu-((ma-b)^q+(ma-b))=0, \\
v^q+(u+m^q)v^{q-1}+u^q+m=0, \\
v+u+m^q-(u+m^q)w^{q-1}=0.
\end{cases}
\end{equation}
$Gal(M|K)\cong \PGL(2,q)$ is generated by the following three automorphisms defined over $\mathbb{F}_{q}$ of order $2,p$ and $q-1$ respectively.
\begin{equation*}
\begin{array}{llll}
\varphi(m)=m, & \varphi(u)=v+u, & \varphi(v)=-v, & \varphi(w)=w^{-1},\\
\varphi(m)=m, & \varphi(u)=u, & \varphi(v)=vw(w+1)^{-1}, & \varphi(w)=w+1,\\
\varphi(m)=m, & \varphi(u)=u, & \varphi(v)=v, & \varphi(w)=\lambda w, \lambda\in \mathbb{F}_{q}^*.\\
\end{array}
\end{equation*}
In particular, $Gal(M|K)$ is defined over $\mathbb{F}_{q^2}$.
\end{theorem}
\begin{proof} We make some computation to show that the above maps are automorphisms of $M$. We begin with the first one.
$$
\begin{array}{lll}
(v+u)^{q+1}+m^q(v+u)^q+m(v+u)-((ma-b)^q+(ma-b))-(u^{q+1}+\\
m^qu^q+mu-((ma-b)^q+(ma-b))=v^{q+1}+v^q(u+m^q)+(u+m^q)v=0.\\
(-v)^q+(u+v+m^q)(-v)^{q-1}+(v+u)^q+m=\\
-(v^q-(u+v+m^q)v^{q-1}-(v^q+u^q+m))=0.\\
-v+(v+u)+m^q-(v+u+m^q)w^{-(q+1)}=\\
(u+m^q)-(v+u+m^q)w^{-(q+1)}=0\\
\end{array}
$$
Now, the computation for the second map.
$$
\begin{array}{lllll}
(vw)^q(w+1)^{-q}+(u+m^q)(vw)^{(q-1)}(w+1)^{-(q-1)}+u^q+m=\\
(w+1)^{-q}((vw)^q+(u+m^q)(vw)^{q-1}(w+1)+(u^q+m)(w^q+1))=\\
(w+1)^{-q}(v^{q-1}w^q(v+u+m^q)+(u+m^q)v^{q-1}w^{q-1}+(u^q+m)w^q+u^q+m)=\\
(w+1)^{-q})((w^q(v^q+(u+m^q)v^{q-1}+u^q+m)+(u+m^q)v^{q-1}w^{q-1}+u^q+m)=\\
(w+1)^{-q}(w^q(v^q+(u+m^q)v^{q-1}+u^q+m)+(u+m^q)+v^{q-1}w^{q-1}+u^q+m)=\\
(w+1)^{-q}((u+m^q)v^{q-1}w^{q-1}+u^q+m)=\\
(w+1)^{-q}(v^{q-1}(v+u+m^q)+u^q+m)=0.\\
(vw)(w+1)^{-1}+u+m^q-(u+m^q)(1+w)^{q-1}= \\
(w+1)^{-1}(vw+(u+m^q)(w+1)-(u+m^q)(1+w)^q)=\\
(w+1)^{-1}(vw+(u+m^q)w-(u+m^q)w^q+u+m^q-(u+m^q)=\\
(w+1)^{-1}(w(v+u+m^q -(u+m^q)w^{q-1})=0.
 \end{array}
$$
Finally, for the third map.
$$
\begin{array}{lll}
v+u+m^q-(u+m^q)(\lambda w)^{q-1}=v+u+m^q-(u+m^q)w^{q-1}=0.
\end{array}
$$
The group generated by the first and the third automorphisms is a dihedral group of order $2(q-1)$ which is a maximal subgroup of $\PGL(2,q)$. Thus, these together with the second automorphism generate the whole $Gal(M|K)=\PGL(2,q)$.
\end{proof}

\begin{lemma}
\label{lemA201222} Let $\cF$ be the affine algebraic curve in $AG(3,\mathbb{F})$ with coordinates $(X,Y,Z)$ defined by
\begin{equation}
\label{eqB20122}
\begin{cases}
F_1=X^{q+1}+Y^qX^q+YX-((Ya-b)^q+(Ya-b))=0, \\
F_2=Z^q+(X+Y^q)Z^{q-1}+X^q+Y=0, \\
\end{cases}
\end{equation}
where $a^{q+1}+b^q+b\neq 0$.
Let $S=(\xi,\eta,\zeta)$ be a point of $\cF$ such that $\xi\notin \mathbb{F}_{q^2}$ and $\zeta \neq 0$. If either $\xi^q+\eta=0$ or $\xi+\eta^q=0$ then $S$ is a non-singular point of $\cF$.
\end{lemma}
\begin{proof} The Jacobian matrix of $\cF$ is
$$
\left(
  \begin{array}{ccc}
    \partial F_1/\partial X & \partial F_1/\partial Y& \partial F_1/\partial Z \\
    \partial F_2/\partial X & \partial F_2/\partial Y& \partial F_2/\partial Z \\
  \end{array}
\right)
=
\left(
  \begin{array}{ccc}
    X^q+Y & X-a & 0 \\
    Z^{q-1} & 1 & -Z^{q-2}(X+Y^q) \\
  \end{array}
\right)
$$
Therefore, if $S$ is singular then $J$ has rank $1$, that is,
\begin{equation}
\label{eqA201222}
\begin{cases}
\xi^q+\eta-(\xi-a)\zeta^{q-1}=0\\
(\xi^q+\eta)(\xi+\eta^q)\zeta^{q-2}=0\\
(\xi-a)(\xi+\eta^q)\zeta^{q-2}=0.
\end{cases}
\end{equation}
We begin with the case $\xi^q+\eta=0$. By $\zeta \ne 0$ the
first equation in (\ref{eqA201222}) yields $\xi=a$ whence $\eta=-a^q$ follows. Therefore, from the first equation in (\ref{eqB20122}), we get $a^{q+1}+b^q+b=0$, a contradiction. Now, $\xi+\eta^q=0$ is assumed.
From the first equation in (\ref{eqB20122}), $\zeta^q+\xi^q+\eta=0$. Since $\zeta \neq 0$, this together with  the second equation in (\ref{eqA201222}) yield $\zeta+\xi-a=0$ whence $\zeta^q+\xi^q-a^q=0$ follows.
Thus, $\eta=-a^q$. Now, the first equation in (\ref{eqB20122}) yields $a^{q+1}+b^q+b=0$, a contradiction.
\end{proof}

\subsubsection{Case  $a,b\in \mathbb{F}$ with $a^{q+1}+b^q+b=0$}
\label{subs2}
For this choice of $a,b$, the polynomial $F(T)$ defined in (\ref{eq1020322}) is reducible as $a^{q+1}+b^q+b=0$ yields $F(a)=0$. Replacing $T$ by $T+a$, $F(T)$ becomes $T(T^q+(a+m^q)T^{q-1}+(m+a^q))$.  Then dividing it by $T$,
we obtain
$$g(T)=T^q+(a+m^q)T^{q-1}+(m+a^q)\in K[T].$$
\begin{lemma}
\label{lemA291122}
$g(T)$ is irreducible over $K$.
\end{lemma}
\begin{proof}
Since $g(T)$ and $h(T)=(m+a^q)T^q+(a+m^q)T+1\in K[T]$ are simultaneously irreducible, it is enough to show that the plane curve $\cU$ of affine equation $U(X,Y)=(Y+a^q)X^q+(a+Y^q)X+1=0$ is non-singular. Since $\deg(\cU)=q+1$ and the line $\ell_\infty$ meets $\cU$ in $q+1$ pairwise distinct points, these points are non-singular. Furthermore, $U_X=Y^q+a$ and $U_Y=X^q$. Since the system consisting of the equations $U(X,Y)=0,Y^q+a=0,X^q=0$ has no solution, the claim follows.
\end{proof}
Let $u$ be a root of $g(T)$ in some overfield of $K$. Comparison of $g(T)$ with (\ref{eqA251122}) shows that the computation in Section \ref{subs1} carried out to obtain $f_2(T)$  from $f_1(T)$ can also be used to find the first
Abhyankar's derivatives $g_1(T)$ and $g_2(T)$ of $g(T)$. From (\ref{eq251122}),
$$\Phi(W)=(u+a+m^q)-(a+m^q)W^{q-1}.$$ Therefore, the first
Abhyankar's derivative $g_1(T)$ of $g(T)$ is
$$g_1(T)=v^{q-2}((u+a+m^q)T^{q-1}-(a+m^q)(T+u)^{q-1}).$$
Moreover, the first Abhyankar's derivative $\Phi_1(T)$ of $\Phi(T)$ at $v$ is,
$$\Phi_1(W)=-\frac{a+m^q}{v}\prod_{i=1}^{q-1}(W-\theta^i v).$$
From this, $\Phi_1(W)$ (and hence $g_2(T)$) is a completely reducible polynomial over $K(u,v)=\mathbb{F}(m,u,v)$.
\begin{theorem}
\label{teoB} Let $a\in \mathbb{F}$. Suppose that $g_1(T)$ is irreducible over $K(u)$. Then $M$ coincides with $K(u,v)=\mathbb{F}(m,u,v)$, and the extension $K(u,v)|K(u)$ is cyclic of degree $(q-1)$. Therefore $M=\mathbb{F}(m,u,v)$
with
\begin{equation}
\label{eq261122A}
\begin{cases}
u^q+(a+m^q)u^{q-1}+a^q+m=0, \\
u+a+m^q-(a+m^q)v^{q-1}=0.
\end{cases}
\end{equation}
\end{theorem}
Moreover, $Gal(M|K)$ is a sharply $2$-transitive group on $\Delta$ such that the $1$-point stabilizer is cyclic.
From Zassenhaus' theorem \cite{Za}, $Gal(M/K)\cong \AGL(1,q)$, and $Gal(M/K)$ acts on $\Delta$ as $\AGL(1,q)$ on the affine line over $\mathbb{F}_q$.

\subsection{The rational BKS case}
\label{ssHon}
We use the same method as for the Hermitian case. For this purpose, it is useful to replace  $(2m-1)/2m$ by $m$. Then (\ref{eq3jul2022C}) becomes the monic polynomial
\begin{equation}
\label{eq3jul2022C1} f=f(T)=T^{q+1}+mT^q+mT-(b-2)(m-1)-\ha a +1\in \mathbb{F}[T].
\end{equation}
\subsubsection{Case  $a,b\in \mathbb{F}$ with $b^{q+1}-(a^q+a)+(b^2-4a)^{(q+1)/2}\neq 0$}
\label{subs1A}
In this case $P(a,b)\not \in \cC$ and hence there exists no $t\in \mathbb{F}$ such that $a=2(t+1)^{q+1}$ and $b=t^q+t+2$.  Therefore, Lemma \ref{lem5dic} applies, and gives the following result.
\begin{lemma}
\label{lemA270412}
$f(T)$ is irreducible over $K$.
\end{lemma}
Let $u$ be a root of $f(T)$ in some extension of $K$ and put $L=K(u)$. A straightforward computation shows that the first Abhyankar's skew derivative of $f$ at $u$ is
\begin{equation}
\label{eqA05dic22}
f_1=f_1(T)=\frac{f(T+u)-f(u)}{T}=T^q+(u+m)T^{q-1}+u^q+m \in L[X].
\end{equation}
Computation to obtain the second Abhyankar's skew derivative of $f(T)$ can be carried out as in Section \ref{subs1} determining first $\Phi(W)$.
\begin{equation}
\label{eq2051222}
\Phi(W)=(v+u+m)-(u+m)W^{q-1}.
\end{equation}
Therefore the second Abhyankar's skew derivative of $f(T)$ at $v$
$$f_2(T)=v^{q-2}((v+u+m)T^{q-1}-(u+m)(T+v)^{q-1}).$$
Furthermore,if $f_1(T)$ is irreducible then the arguments on the third Abhyankar's skew derivative used in Section \ref{subs1} remain valid in the present case  whenever $m^q$ is replaced by $m$. In particular,
$f_3(T)$ is completely reducible over $K(u,v)=\mathbb{F}(m,u,v)$. Therefore, the following result holds.
\begin{theorem}
\label{teoC} Suppose that $f_1(T)$ is irreducible over $K=\mathbb{F}(m)$. Then $f_2(T)$ is irreducible over $K(u)$, $M$ coincides with $K(u,v,w)=\mathbb{F}(m,u,v,w)$ where the extension $K(u,v,w)|K(u,v)$ is cyclic of degree $q-1$, and $M=\mathbb{F}(m,u,v,w)$ with
\begin{equation}
\label{eqC051222}
\begin{cases}
u^{q+1}+mu^q+mu-(b-2)(m-1)-\ha a+1=0, \\
v^q+(u+m)v^{q-1}+u^q+m=0, \\
v+u+m-(u+m)w^{q-1}=0.
\end{cases}
\end{equation}
\end{theorem}
\subsubsection{Case  $a,b\in \mathbb{F}$ with $b^{q+1}-(a^q+a)+(b^2-4a)^{(q+1)/2}= 0$}
\label{ssBKSon}
In this case, $P(a,b)\in \cC$ and hence  there exists $t\in \mathbb{F}$ such that $a=2(t+1)^{q+1}$ and $b=t^q+t+2$. Replace $T$ by $T+t$ in (\ref{eq3jul2022C1}). Then $f(T)=Tg(T)$ where
\begin{equation}
\label{eqA11122022Y}
g(T)=T^q+(m+t)T^{q-1}+m+t^q.
\end{equation}
From now on we assume $P(a,b)\notin \PG(2,q)$, i.e. $t\notin \mathbb{F}_{q^2}$. From Lemma \ref{lem7dic} the following claim follows.
\begin{lemma}
\label{lem11122022A}$g(T)$ is irreducible over $K$.
\end{lemma}
Now, the arguments in Section \ref{subs2} remain valid if $a$ is replaced by $t$. Therefore,
the following result is obtained.
\begin{theorem}
\label{teoD} Let $t\in \mathbb{F}\setminus \mathbb{F}_{q^2}$. Suppose that $g_1(T)$ is irreducible over $K(u)$. Then $M$ coincides with $K(u,v)=\mathbb{F}(m,u,v)$, and the extension $K(u,v)|K(u)$ is cyclic of degree $(q-1)$. Therefore $M=\mathbb{F}(m,u,v)$
with
\begin{equation}
\label{eq261122D}
\begin{cases}
u^q+(t+m)u^{q-1}+t^q+m=0, \\
u+t+m-(t+m)v^{q-1}=0.
\end{cases}
\end{equation}
\end{theorem}
Moreover, $Gal(M|K)$ is a sharply $2$-transitive group on $\Delta$ such that the $1$-point stabilizer is cyclic.
From Zassenhaus' theorem \cite{Za}, $Gal(M/K)\cong \AGL(1,q)$, and $Gal(M/K)$ acts on $\Delta$ as $\AGL(1,q)$ on the affine line over $\mathbb{F}_q$.

In the theorems of this section, Theorems \ref{teoA}, \ref{teoB}, \ref{teoC} and \ref{teoD}, it is assumed that the first Abhyankar derivative is irreducible. Actually, this hypothesis can be dropped. Our proof requires a further tool, namely a classical theorem due to van der Waerden, and it is detailed in the following section.
\section{van der Waerden's theorem}
\label{vWt}
In this section,  $L|K$ stands for a finite separable extension of function fields, $M$ for its splitting field (equivalently, for its Galois closure), and $G= \Gal(M|K)$ for the Galois group. Also, $A=\Gal(M|L)$ denotes the Galois group of the Galois extension $M|L$. Furthermore, if $P$ is a place of $K$ and $\cS$ is a the set of all places of $L$ over $P$, then $e(S|P)$ and $f(S|P)$ denote the ramification index and the relative degree for $S\in \cS$, and $W$ is a place of $M$ lying above $P$. Moreover, $D(W|P)$ is the decomposition group and $I(W|P)$ is the inertia group. Finally, $f(x)$ denotes an irreducible polynomial over $K$ such that $L=K(\alpha_1)$ for $f(\alpha_1)=0$.

%Since $M$ is the splitting field of $f(x)$ and $f(x)$ is separable, if $\{\alpha_1,\ldots,\alpha_n\}$ is the set of all roots of $f(X)$ in $M$ then $\deg(f(x))=[L:K]$ and $|G|=n|A|$ as the roots of $f(x)$ are pairwise distinct and have multiplicity one. Moreover, $G$ has a faithful action on $\{\alpha_1,\ldots,\alpha_n\}$ which is transitive, that is, $G$ can be viewed a permutation group on $\{\alpha_1,\ldots,\alpha_n\}$.

Artin and later on van der Waerden investigated the action of $D(W|P)$ and $I(W|P)$ on the set of all roots of $f(x)$ and, in particular, how the their orbits are linked to the ramification picture of $\cS$. The following result is due to van der Waerden.
\begin{result}
\label{vWsatz} \rm{(\cite[Satz I]{vW})} {\em{ Under the action of $D(W|P)$, the set of all roots of $f(x)$ splits into as many as $|\cS|$ orbits. Each $D(W|P)$-orbit consists of $e(S|P)f(S|P)$ roots of $f(x)$ while each $I(W|P)$-orbit does of $e(S|P)$. }}
\end{result}
In particular, if $L|K$ is ramified at $W$ then $D(W|P)$ is non-trivial.

We reproduce the proof of Result \ref{vWsatz} using  notation and terminology from \cite{HKT,sti}.

Let $S\in \cS$, and take a place $U$ of $M$ lying over $S$, Since $S$ lies over $P$, we have that $U$ is in the $G$-orbit of $W$. Wherefore, there exists a $g_S\in G$ (not uniquely determined) such that $g_S(P)=U$. Clearly, the cosets $Ag_S$ with $S$ running over $\cS$ give a partition of the $G$-orbit of $W$. Now, consider the sets $Ag_SD(W|P)$ with $S$ ranging over $\cS$.

We show that if $g\in Ag_SD(W|P)$ for some $S\in \cS$ then $g^{-1}(\alpha_1)$ and $g_S^{-1}(\alpha_1)$ fall in the same $D(W|P)$-orbit on $\{\alpha_1,\ldots,\alpha_n\}$. Write $g=ag_Sd$ with $a\in A$ and  $d\in D(W|P)$. When $g^{-1}=
d^{-1}g_S^{-1}a^{-1}$, and hence $g^{-1}(\alpha_1)=d^{-1}(g_S^{-1}(a^{-1}(\alpha_1)))=d^{-1}(g_S^{-1}(\alpha_1))$. Therefore, $g^{-1}(\alpha_1)$ can be viewed as the image of $g_S^{-1}(\alpha_1)$ by $d^{-1}$.  The converse also holds. In fact, let $g^{-1}(\alpha_1)=d^{-1}(g_S^{-1}(\alpha_1))$. Then $(gd^{-1}g_S^{-1})(\alpha_1)=\alpha_1$. Therefore, $gd^{-1}g_S^{-1}$ fixes $L$ element-wise as $L=K(\alpha_1)$. Thus $gd^{-1}g_S^{-1}\in A$ whence
$g=ag_Sd\in Ag_SD(W|P)$.

Since $G$ is transitive on $\{\alpha_1,\ldots,\alpha_n\}$, each  $\alpha_i=g(\alpha_1)$ for some $g\in G$. Therefore there exists a bijective correspondence between the places $S\in \cQ$ and the $D(W|P)$-orbits on
$\{\alpha_1,\ldots,\alpha_n\}$. Since $D(W|P)$ is not transitive in general on $\{\alpha_1,\ldots,\alpha_n\}$ one can ask to compute the size of such a $D(W|P)$-orbit on $\{\alpha_1,\ldots,\alpha_n\}$.
%$$\frac{|D(R|P)|}{|D(R|P)_{\alpha}|$$

 We show that $e(S|P)f(S|P)$ is the size of the corresponding $D(W|P)$-orbit on $\{\alpha_1,\ldots,\alpha_n\}$. Wake $\alpha\in \{\alpha_1,\ldots,\alpha_n\}$. Let $D(W|P)_\alpha$ be the stabilizer of $\alpha$ in $D(W|P)$.
 Then the size of the $D(W|P)$-orbit of $\alpha$ is given by $|D(W|P)|/|D(W|P)_\alpha|$. We may think about $D(W|P)_{\alpha}$ as the intersection of $D(W|P)$ with $G_\alpha$. Take $g\in G$ such that $g^{-1}(\alpha_1)=\alpha$. Then
 $g^{-1}Ag(\alpha)=\alpha$ and hence $g^{-1}Ag$ fixes $K(\alpha)$ element-wise. Wherefore, $G_\alpha=g^{-1}Ag$. This yields
 $$ \frac{|D(W|P)|}{|D(W|P)_{\alpha}|}=\frac{|D(W|P)|}{|D(W|P)\cap g^{-1}Ag|}$$
 Let $g(W)=U$. When $D(W|P)=g^{-1}D(U|P)g$, and hence $D(W|P)\cap g^{-1}Ag=g^{-1}(D(U|P)\cap A) g$. Since $|D(W|P|)=D(U|P)$, we have
 $$ \frac{|D(W|P)|}{|D(W|P)_{\alpha}|}=\frac{|D(U|P)|}{|D(U|P)\cap A|}.$$
 Since $M|K$ is a Galois extension, $|D(U|P)|=e(U|P)f(U|P)$. Moreover, $D(U|P)\cap A$ is the stabilizer of $A$ at $U$. Since $M|L$ is a Galois extension, we also have $D(U|P)\cap A=e(U|S)f(U|S)$. Therefore,
 the $D(W|P)$ orbit of $\alpha$ has size $e(W|P)f(W|P)/e(W|Q)f(W|Q)=e(S|P)f(S|P)$ whence the first claim follows.

Each $D(S|P)$-orbits on $\{\alpha_1,\ldots,\alpha_n\}$ splits further into $I(S|P)$-orbits. The above argument applied to $I(W|P)$ gives that the size of $\alpha$ under the action of $I(D|P)$-orbit is
$$ \frac{|I(W|P)|}{|I(W|P)_{\alpha}|}=\frac{|I(U|P)|}{|I(U|P)\cap A|}=\frac{e(U|P)}{e(U|Q)}=e(Q|P).$$
This ends the proof of Result \ref{vWsatz}.

Now we look inside the case where $K=\bar{\mathbb{F}}_q(t)$. Then the decomposition and inertia groups coincide. Let $f(x)=x^n+\ldots+a_{n-1}(t)x+a_n(t)$ with $a_i(t)\in \mathbb{F}_q[t]$ and $\deg\, a_i(t)\le i$.  Let $\cC$ be the irreducible (possible singular) plane curve of equation $f(X,T)=X^n+\ldots+a_{n-1}(T)X+a_n(T)=0$. Observe that $X_\infty \not\in \cC$.

The places of $K$ are the points of the projective line $\PG(1,\bar{\mathbb{F}}_q)$. Let $x_1,\ldots,x_{u(\tau)}$ be the roots of the polynomial $f(X,\tau)$ in the indeterminate $X$. Geometrically speaking, let $Q(\tau,x_1),\ldots Q(\tau,x_v(\tau))$ be the centers of the branches $\gamma_1,\ldots,\gamma_{u(\tau)}$ of $\cC$ whose centers lie over $\tau$.
 Note that $v(\tau)\ge u(\tau)$.
 Let $\ell_\tau$ be the vertical line through $\tau$. From B\'ezout's theorem, $$\sum_{i=1}^{v(\tau)} I(\gamma_i\cap \ell_\tau,Q_i(\tau))=\deg(\cC).$$
Therefore, $\sum_{i=1}^{v(\tau)} e(\gamma_i|\tau)=\deg(f(X,T))$ since $I(\gamma_i\cap \ell_\tau,Q_i(\tau))$ is the ramification index of $e(\gamma_i|\tau)$.
Result \ref{vWsatz} shows that the $D(R|\tau)$ orbits on the set $\{\alpha_1,\ldots,\alpha_n\}$ are as many as the $v(\tau)$ orbits and that their sizes are $ord (\gamma_1),\ldots,ord (\gamma_v(\tau))$. In particular,
if $\gamma_i$ is linear then $D(P|\tau)$ has a fixed point. This occurs in particular when $Q_i(\tau,\alpha_i)$ is a non-singular point and the tangent to $\cC$ at that point is not the vertical line. It turns out
that if each of the common points of $\cC$ and $\ell$ is non-singular and $\ell$ is tangent to $\cC$ at exactly one point then $D(R|P)$ is a transposition on $\{\alpha_1,\ldots,\alpha_n\}$.

If we drop the hypothesis  $\deg\,a_i(t)\le i$, some changes are needed. Assume that $X_\infty$ is a point of $\cC$. Let $\gamma$ be a branch of $\cC$ centered at $Y_\infty$ with a primitive representation $(x=x(t), y=y(t))$, Two cases arise according as $\gamma$ is a pole or not of $x$. In the latter case $x(t)=c_it^i+c_jt^j+\ldots$ with $i\ge 0$ (and $ord (y(t))<0$). If $i=0$ then $x=c_0$ is the tangent of $\gamma$, and  $\gamma$ lies over the place $c_0$ of $K$, and we have $e(\gamma|c_0)=j$. If $i\ge 1$ then $x=0$ is the tangent of $\gamma$, and $\gamma$ lies over the place $0$ of $K$ and $e(\gamma|0=i$. Otherwise, $i<0$ and the tangent line of $\gamma$ is the line at infinity and $\gamma$ lies over the place $\infty$ of $K$, and we have $e(\gamma|\infty)=-i$.

\subsection{Example, Hermitian curve I}
Let $\cC$ be the Hermitian curve with affine equation $X^{q+1}+T^{q+1}+1=0$. Let $K=\bar{F}_q(t)$ and $L=K(x)$ with $x^{q+1}+t^{q+1}+1=0$. Then $L$ is a Galois extension of $K$ (i.e. $M=L$) as $Gal(L|K)$ consists of all automorphisms $\alpha$ of $\cC$ with $\alpha(t)=t,\alpha(x)=\lambda x$ and $\lambda^{q+1}=1$. In particular, $Gal(L|K)$ acts on the set of all roots of the polynomial $f(X)=X^{q+1}+t^{q+1}+1\in K[X]$ as a sharply transitive permutation group. Let $t_0\in \bar{F}_q(t)\cup \{\infty\}$. Then vertical line $\ell$ through the point $U=(1:0:0)$ meets $\cC$ in pairwise distinct points except for $t_0^{q+1}=0$ in which case $\ell$ is the tangent to $\cC$ and $I(\cC\cap\ell,U)=q+1$ so that $\cC\cap \ell=\{U\}$. From this van der Waerden's theorem follows, since the stabilizer of any place of $M$ in $Gal(L|K)$ is trivial.
\subsection{Example, Hermitian curve II}
Let $\cC$ be the Hermitian curve with affine equation $X^q-X-\omega T^{q+1}=0$ with $\omega^{q-1}=-1$. Let $K=\bar{F}_q(t)$ and $L=K(x)$ with $x^q-x-\omega t^{q+1}=0$. Then $L$ is a Galois extension of $K$ (i.e. $M=L$) as $Gal(L|K)$ consists of all automorphisms $\alpha$ of $\cC$ with $\alpha(t)=t,\alpha(x)=x+a$ and $a\in \mathbb{F}_q$. In particular, $Gal(L|K)$ acts on the set of all roots of the polynomial $f(X)=X^q-X-t^{q+1} \in K[X]$ as a sharply transitive permutation group. Let $t_0\in \bar{F}_q(t)$. Then vertical line $\ell$ through the point $U=(t_0,0)$ meets $\cC$ in $q$ pairwise distinct points. $\cC$ has a just one point at infinity, namely $X_\infty$. The unique branch $\gamma_\infty$ of $\cC$ centered at $X_\infty$ has a primitive representation $(t=s^{-q},s^{-1}+\ldots)$ and hence $e(\gamma_\infty|\infty)=q$. This yields van der Waerden's theorem, as the stabilizer of any place of $M$ in $Gal(L|K)$ is trivial.
 \subsection{Example, Hermitian curve III}
  \label{ssHIII} Let $\cC$ be the curve with affine equation $X^q+T^qX+T=0$. It is known that $\cC$ is (linearly) isomorphic to the Hermitian curve over $\mathbb{F}_{q^6}$. Let $K=\bar{\mathbb{F}}_q(t)$ and $L=K(x)$ with $x^q+t^qx+t=0$. Then the automorphism group of $\cC$ fixing the points $(t_0,0)$ of $\cC$ is trivial. Therefore, $L|K$ is not a Galois extension. Moreover, $\cC$ has exactly two points at infinity, namely $X_\infty$ and $T_\infty$. The unique branch $\gamma_\infty$ of $\cC$ centered at $T_\infty$ has a primitive representation $(t=s^{-1},x=-s^{(q-1)}(1+\ldots))$ while the unique branch $\delta_\infty$ of $\cC$ centered at $X_\infty$ has a primitive representation $(t=s^{-(q-1)}(1+\ldots),x=s^{-q}(1+\ldots))$. Therefore, $L$ has exactly two places (branches) lying over $\infty$ and $e(\gamma_\infty|\infty)=1$, $e(\delta_\infty|\infty)=q-1$. Van der Waerden's theorem shows that $Gal(M|K)$ has a subgroup that has two orbits on the set $\{\alpha_1,\ldots,\alpha_q\}$ of all roots of $f(X)=X^q+t^qX+t=0\in K[x,]$ of size $1$ and $q-1$ respectively. Since $Gal(M|K)$ is
transitive on $\{\alpha_1,\ldots,\alpha_q\}$ this yields that $Gal(M|K)$ is doubly transitive on  $\{\alpha_1,\ldots,\alpha_q\}$. We show that $Gal(M|K)$ is sharply $2$-transitive on $\{\alpha_1,\ldots,\alpha_q\}$, and hence
$Gal(M|K)$ is the semidirect product of an elementary abelian group of order $q$ by a cyclic complement of order $q-1$. First we give an explicit equation for $M$. Let $N$ be an extension of $L$ of degree $q-1$ defined by  $\mathbb{K}(N)=\mathbb{K}(x,t,y)$ with $x^q+xt^q+t=0$, and $y^{q-1}=t^q$. Clearly, the map $\varphi_\lambda:(x,t,y)\mapsto (x,t,\lambda y)$ with $\lambda\in {\mathbb{F}_q}^*$ is an automorphism of $N$ which fixes $L$ element-wise, and $\Lambda=\{\varphi_\lambda|\lambda\in {\mathbb{F}_q}^*\}$ is (cyclic) group of order $q-1$. Since $[N:L]=q-1$, this shows that $[N|L$ is a Galois extension with $Gal(M|L)=\Lambda$. Furthermore, for a fixed
$\alpha\in \mathbb{F}_{q^2}$ with $\alpha^{q-1}=-1$, the map $\sigma_\alpha: (x,t,y)\mapsto (x+\alpha y,t,y)$  is an automorphism of $N$ which fixes $t$, and hence $K$ element-wise. In fact, $(x+\alpha y)^q+(x+\alpha y)t^q+t=x^q+xt^q+t+\alpha y(\alpha^{q-1}y^{q-1}+t^q)=0$. Therefore, $\Sigma=\{\sigma_\alpha|\alpha^{q-1}=-1\}\cup \{id\}$ is an elementary abelian group of order $q$. Moreover, $\Sigma$ together with $\Lambda$ generate a group $H$ of order $q(q-1)$ which is the semidirect product of $\Sigma$ with complement $\Lambda$. Clearly, $H$ fixes $t$, and hence $K$ element-wise. Let $R$ be the fixed field of $H$, Then $q(q-1)=|H|=[N:R]$. Since $[N:K]=q(q-1)$, this yields $K=R$, that is, $K$ is the fixed field of $H$. As $K\le L \le N$, this shows that $M\le N$ up to a birational isomorphism. On the other hand, $|Gal(M|K)|\ge q(q-1)$, as $Gal(M|K)$ is a 2-transitive permutation group of degree $q$, Therefore $N=M$, up to a birational isomorphism.  Eliminating $t$ from the equations defining $N$ shows that $N=\mathbb{K}(x,y)$ with $y^{q-1}+x^{q^2}+x^qy^{q(q-1)}=0$. Replacing $xy^{-1}$ by $\xi$ and $y^{-1}$ by $\eta$ shows that $M$ is birationally isomorphic to $\mathbb{K}(\xi,\eta)$ with $\eta^{q^2-q+1}+\xi^{q^2}+\xi^q$, and hence $M\cong \mathbb{K}(\xi,\eta); \eta^{q^2-q+1}+\xi^q+\xi=0.$  By a result of Tafazolian and Torres \cite{TT}, $M$ is $\mathbb{F}_{q^6}$-covered by the Hermitian function field $\mathbb{F}_{q^6}(\cH_{q^3})$. In particular, $M$ is $\mathbb{F}_{q^6}$ maximal of genus $g=\ha q(q-1)^2$ and is Galois subcover of  $\mathbb{F}_{q^6}(\cH_{q^3})$.
%Our next topic is motivated by a question raised by Artin and van der Waerden in 1920s.
\section{Case $P\notin \Omega$}
\label{ssnotin}
We keep up our notation from Sections \ref{abhy} and \ref{vWt}.
\subsection{Hermitian case}
\label{vWvher}
We show that the irreducibility condition in Theorem \ref{teoA} is fulfilled.
%\subsection{Case  $a,b\in \mathbb{F}$ with $a^{q+1}+b^q+b\neq 0$}
We first prove that $Gal(f,K)$ is $2$-transitive on $\Delta$. From Lemma \ref{lemA271122}, $Gal(f,K)$ is transitive on $\Delta$. Therefore, it is sufficient to prove that the $1$-point stabilizer of $Gal(f,K)$ on $\Delta$ is transitive on the remaining $q$ points.

Consider the function field $L=\mathbb{F}(m,u)$ with $u^{q+1}+m^qu^q+mu-(ma-b)^q+(ma-b))=(u^q+m)(u+m^q)-(m^{q+1}+(ma-b)^q+ma-b)=0$, as the algebraic extension $L|K$ of degree $q+1$. Let $m_i\in \mathbb{F}$ with $i=1,2,\ldots,q+1$ be the roots of the polynomial $Y^{q+1}+(Ya-b)^q+Ya-b\in \mathbb{F}[Y]$. The results stated in Section \ref{hhc} show that the tangents to $\cH_q$ passing through the point $Q(a,b)$ are exactly the lines $\ell_i$ of equation $Y=m_i(x-a)+b$. Furthermore, the tangency point on $\ell_i$ is $P_i=P_i(-\sqrt[q]{m_i},m_i(\xi-a)+b)$ with $I(P_i, \cH_q\cap \ell_i)=q$, and  the remaining intersection of $\ell_i$ with $\cH_q$ is the point $R_i=(-m_i^q,m_i(-m_i^q-a)+b)$ with $I(R_i,\cH_q\cap \ell_i)=1$. Let $\cS_i$ be the set of places of $K(u)$ lying over $m_i$ in the covering $K(u)|K$. Then $\cS_i=\{P_i,R_i\}$ and $e(P_i|m_i)=q$ and $e(R_i|m_i)=1$. Let $W_i$ be a place of $M$ lying over $m_i$ in the covering $M|K$. From Result \ref{vWsatz}, the inertia group $I(W_i|m_i)$ has two orbits on the set $\Delta$ of the roots of the polynomial (\ref{eq1020322}), one of size $q$ and another of size $1$. Therefore, $I(T_i|m_i)$ as a subgroup of $Gal(f,K)$ acting on $\Delta$ fixes a point and transitive on the remaining points.

By the results recalled in Section \ref{subs1}, the $2$-transitivity of $Gal(f,K)$ has the following implication.
\begin{proposition}
\label{prop271122} The first Abhyankar's skew derivative $f_1(T)$ of $f(T)$ given in (\ref{eq1020322}) is an irreducible polynomial over $K(u)$, and the irreducibility condition in Theorem \ref{teoA} can be dropped.
\end{proposition}
Therefore, Theorem \ref{teoA} applies and it determines the structure of $Gal(M|K)$. In fact, $Gal(M|K)$ turns out be a sharply $3$-transitive group on $\Delta$ such that the $2$-point stabilizer is cyclic.
From Zassenhaus' theorem \cite{Za}, $Gal(M/K)\cong \PGL(2,q)$, and $Gal(M/K)$ acts on $\Delta$ as $\PGL(2,q)$ on the projective line over $\mathbb{F}_q$, in its unique $3$-transitive permutation representation.

We go on by describing the set $\cW_i$ of places of $M$ lying over the place $m_i$ of $K=\mathbb{F}(m)$. As we have observed, the set $\cS_i$ of the places of $K(u)$ lying over $(m_i)$ consists of two points $P_i$ and $R_i$ unless $m_i\in \mathbb{F}_{q^2}$ in which case $P_i=Q_i$ and hence $S_i$ is reduced in a unique place. We treat these two cases separately using the equations in (\ref{eq261122})

\subsubsection{Case $m_i\notin \mathbb{F}_{q^2}$.} We begin with $P_i$. The second equation in (\ref{eq261122}) together with Lemma \ref{lemA201222} show that there exist exactly two places in $K(u,v)$ lying over $P_i$, namely those with center at $P_{i,1}=(m_i,u_i,0)$ and $P_{i,2}$ by $(m_i,u_i,-u_i-m^q)$ respectively. Here $e(P_{i,1}|P_i)=q-1$ and $e(P_{i,2}|P_i)=1$. From the third equation in (\ref{eq261122}), there exist exactly $q-1$ pairwise distinct places of $M$ lying over $P_{i,1}$, they are centered at $P_{1,i,j}=(m_i,u_i,0,w_i^j)$ where $w_i\in \mathbb{F}_q$ is a fixed primitive ($q-1$)-st root of unity. Here $e(P_{i,1,j}|P_{i,1})=1$. Let $\kappa$ be the number of places of $M$ lying over $P_{i,1}$. Then $\le \kappa \le q-1$, and all these places are centered at $(m_i,u_i,-u_i-m_i^q,0)$. We repeat the same argument for $R_i$ which is the place of $K(u)$ centered at $(m_i,u_i)$ with $u_i+m_i^q=0$. From the second equation in (\ref{eq261122}), there exists a unique place of $K(u,v)$ lying over $R_i$, identified by $R_{i,1}=(m_i,u_i,-\sqrt[q]{u_i^q+m_i})$ where $e(R_{i,1}|R_i)=q$.
Let $\rho$ be the number of places of $M$ lying over $R_{i,1}$. Then $\le\rho\le q-1$ and all these places are centered at the same point. Therefore, there are as many as $q-1+\kappa+\rho$ places of $M$ lying over $(m_i)$. They form $\cW_i$ which is an orbit under the action of $Gal(M|K)$. Since $Gal(M|K)\cong \PGL(2,q)$ and $1$-point stabilizer of $Gal(M|K)$ contains a subgroup of order $q$, the Dickson classification of subgroups of $\PGL(2,q)$ yields that $1$-point stabilizer has order $q(q-1)/r$ with $r\mid(q-1)$. Therefore $|\cW_i|=(q+1)r$ whence $(q+1)r=q-1+\kappa+\rho$ follows. From this $r\leq 2$. For $r=2$, the $1$-point stabilizer of any subgroup of order $q-1$ has order $\ha(q-1)$. On the other hand, case $r=2$ may only occur when $\kappa+\rho=q+3$ and hence one of $\kappa$ and $\rho$ is equal to $q-1$ and the other is $4$. Therefore, if $r=2$ then the subgroup of $Gal(M|K)$ of order $q-1$, namely $Gal(M|K(u,v)$, has an orbit of size $4$ and hence its subgroup of order $\ha(q-1)$ cannot fix a place in $\cW_i$. Thus $r=1$ and $|\cW_i|=q+1$. In particular, $\kappa=\rho=1$.
 \subsubsection{Case $m_i\in \mathbb{F}_{q^2}$.} This time, there exists just one place in $K(u,v)$ lying over $P_i=Q_i$. Let $T$ be a place of $M$ lying over $P_i$ in the covering $M|K$. From Result \ref{vWsatz}, the stabilizer of $T$ in $Gal(M|K)\cong \PGL(2,q)$ is transitive on the set $\Delta$ of the $q+1$ roots of the polynomial (\ref{eq2411A22}), and hence its order is a multiply $r(q+1)$ of $q+1$. From the Dickson classification of subgroups of $\PGL(2,q)$, either $r=1$, or the stabilizer of $T$ contains $PSL(2,q)$. The latter case cannot actually occur, since the $p$-subgroup of an automorphism group of any function field fixing a place is a normal subgroup. Thus the stabilizer of $T$ in $Gal(M|K)$ is a cyclic group of order $q+1$, and hence $|\cW_i|=q(q-1)$.

% Actually, $\cW_i$ for $i=1,\ldots q+1$  are the only short orbits of $Gal(M|K)$, i.e. the covering $M|K$ only ramifies at the points in $\cW=\cW_1\cup \cdots \cup \cW_{q+1}$.

We are in a position to compute the genus $\mathfrak{g}(M)$. For $T\in \cW$, let $G_T^{(k)}$ be the $k$-th ramification group of $G$ at $T$. Since $G_T$ is isomorphic to the $1$-point stabilizer of $\PGL(2,q)$ in its action on the projective line, we have that $G_T$  is a semidirect product of an elementary abelian normal subgroup $S_q$ of order $q$ by a cyclic complement of order $q-1$. Since the non-trivial elements in $S_q$ form a unique conjugacy class in $G_T$, the non-trivial ramification groups with $k\ge 1$ coincide with $S_q$. Thus, for some positive integer $s=s_i$, $i=1,2,\ldots q+1$ depending only on $\cW_i$, $S_q=G_T^{(1)}=\cdots =G_T^{(s_i)}$ and $G_T^{(k)}=\{1\}$ with $k>s_i$. Two cases are treated separately according as there is a tangent to $\cH_q$ at a point in $\PG(2,q^2)$ which passes through $P(a,b)$.
\subsubsection{$m_i\in \mathbb{F}\setminus \mathbb{F}_{q^2}$ for $i=1,2,\ldots,q+1$}
 From the Hurwitz genus formula \cite[Theorem 11.72]{HKT} applied to the Galois covering $M|K$,
\begin{equation}
\label{eqA291122}
\begin{cases}
2\mathfrak{g}(M)-2&= -2|\PGL(2,q)|+\sum_{T\in \cW}\sum_{i\ge 0}(|G_T^{(i)}|-1)\\
&=-2(q^3-q)+(q+1)^2(q(q-1)-1)+(q+1)\sum_{i=1}^{q+1}s_i(q-1))\\
&=q^4-q^3-2q^2-q-1+(q+1)^2(q-1)\sum_{i=1}^{q+1}s_i.
\end{cases}
\end{equation}
The subfield $K(u)=H_q$ of $M$ is a Galois cover of $M$ with Galois group $Gal(M|K(u)$ isomorphic to $\AGL(1,q)$. Since $\mathfrak{g}(H_q)=\ha q(q-1)$,
the Hurwitz genus formula \cite[Theorem 11.72]{HKT} applied to the Galois covering $M|K(u)$ gives
\begin{equation}
\label{eqB291122}
\begin{cases}
2\mathfrak{g}(M)-2&= |\AGL(1,q)|(q(q-1)-2)+\sum_{T\in \cW}\sum_{i\ge 0}(|G_T^{(i)}|-1)\\
&=q(q-1)(q(q-1)-2)+(q+1)((q(q-1)-1)+\\
&\qquad \sum_{i=1}^{q+1}s_i(q-1)+(q+1)q(q-2)).
\end{cases}
\end{equation}
Comparison of (\ref{eqA291122}) with (\ref{eqB291122}) yields $\sum_{i=1}^{q+1}s_i=q+1$. Since $s_i\ge 1$, this yields $s_1=\ldots=s_{q+1}=1$. Therefore,
\begin{equation}
\label{eq291122}
2\mathfrak{g}(M)-2=(q+1)(q^3-q-2)=q^4-q^2-2q-2.
\end{equation}
\subsubsection{$m_i\in \mathbb{F}\setminus \mathbb{F}_{q^2}$ for $i=1,2,\ldots,q$ and $m_{q+1}\in \mathbb{F}_{q^2}$}

\begin{theorem}
\label{teoE} Let $k=q^{2r}+1\pm q^{r+1}(q-1)$ where $\pm$ is taken according as $r$ is even or odd. In $\PG(2,q^{2r})$ with $r\ge 3$, let $\Omega$ be the $(k,q+1)$-arc consisting of all points of the Hermitian curve. If $r>3$ then $\Omega$ is complete.
\end{theorem}
\begin{proof} We show that some long orbit of $Gal(M|K)$ consists of places defined over $\mathbb{F}_{q^{2r}}$. Since $Gal(M|K)$ has exactly $(q+1)$ short orbits, each of size $q+1$, $M$ has as many as $(q+1)^2$ ramified places. By (\ref{eq291122}) the Hasse-Weil lower bound ensures at least $q^{2r}+1-q^r(q^4-q^2-2q)=q^{2r}-q^{r+4}+q^{r+2}+2q^{r+1}$ places of $M$ defined over $\mathbb{F}_{q^{2r}}$. For $r>3$, this number is larger than
$(q+1)^2$. Therefore, as long as $r>3$, $M$ has a unramified place $P_0$ over $\mathbb{F}_{2r}$. Since $Gal(M|K)$ is defined over $\mathbb{F}_{q^{2r}}$ (Theorem \ref{teoA} and Proposition \ref{prop271122}), the long orbit of $P_0$ under the action of $Gal(M|K)$ is entirely consists of places defined over $\mathbb{F}_{q^{2r}}$. Therefore, the place $m_0$ of $K$ lying under $P_0$ has the required property, that is, the roots of the polynomial $F(T)=T^{q+1}+m_0^qT^q+m_0T-((m_0a-b)^q+m_0a-b)$
are pairwise distinct and belong to $\mathbb{F}_{q^{2r}}$.
  \end{proof}
\subsection{The rational BKS case}
Our goal is to show that Proposition \ref{prop271122} holds true for the rational BKS curve unless $Q\in \PG(2,q)$.  For this purpose, we proceed as in Section \ref{vWvher}. First we prove the transitivity of the $1$-point stabilizer of $Gal(f,K)$ on the remaining $q$ points of $\Delta$.

This time, the function field $L$ to be investigated is
$$L=\mathbb{F}(m,u);\, 2mu^{q+1}+(2m-1)u^q+(2m-1)u+m(2-a)+b-2.$$
  The algebraic extension $L|K$ of degree $q+1$. For $a \neq 0$, let $m_i\in \mathbb{F}$ with $i=1,2$ be the roots of the polynomial $2aY^2-2bY+1\in \mathbb{F}[Y]$, i.e. $m_i=(b\pm \sqrt{b^2-2a})/(2a)$, and let $\ell_i$ be the line of equation $Y=m_i(X-a)+b$. For $a=0$, let $\ell_1$ be the line of equation $Y=\ha b^{-1}X+b$ and $\ell_2$ be the line of equation $X=0$.
The results stated in Section \ref{hhc} show that the tangents to $\cC$ passing through the point $Q(a,b)$ are exactly the lines $\ell_i$. Furthermore, $\ell_i$ has exactly two common points with $\cC$, namely the points $P_i,Q_i$ where $I(P_i,\cC\cap \ell_i)=1$ and  $I(Q_i,\cC\cap \ell_i)=q$, unless $P_i=Q_i$ and $I(Q_i,\cC\cap \ell_i)=q+1$. In the latter exceptional case, $m_i\in \mathbb{F}_q$. If this occurs for $i=1,2$ then $Q(a,b)\in \PG(2,q)\setminus \cC$, i.e. $Q(a,b)$ is an external point to $\cC^2$ in $\PG(2,q)$. Dismissing this case allows us to assume that $I(P_1,\cC\cap \ell_1)=1$ and  $I(Q_1,\cC\cap \ell_1)=q$. Therefore we may argue as in Section \ref{vWvher} and prove that the $1$-point stabilizer of $Gal(f,K)$ is transitive on the remaining $q$ points of $\Delta$. From this, $Gal(f,K)$ is a 2-transitive permutation group on $\Delta$. More precisely, $Gal(f,K)\cong \PGL(2,q)$ and it acts on $\Delta$ as $\PGL(2,q)$ in its unique $3$-transitive permutation representation on the projective line over $\mathbb{F}_q$.

Therefore, the following result is obtained.
\begin{proposition}
\label{prop271122BKS}  Assume that $Q\notin \PG(2,q)$. Then the first Abhyankar's skew derivative $f_1(T)$ of $f(T)$ given in (\ref{eqA05dic22}) is an irreducible polynomial over $K(u)$, and the irreducibility condition in Theorem \ref{teoA} can be dropped.
\end{proposition}
For $q=11, r=3$, Magma computation shows that the dismissed case, $Q\in \PGL(2,q)$ is a real exception for Proposition \ref{prop271122BKS}, as  $Gal(f,K)$ has order $2(q+1)$ in that case.

From now on we assume $Q\notin \PG(2,q)$. For $i=1,2$, let $\cW_i$ be the set of places of $M$ lying over the place $m_i$ of $K=\mathbb{F}(m)\cup \{\infty\}$. We distinguish two cases, called general and special, according as $P_i$ and $Q_i$ are distinct or coinciding.

For the general case, we may proceed as in Section  \ref{vWvher}.

Assume first that $m_1\neq m_2$. In this case, $Gal(M|K)$ has exactly two short orbits on $M$, namely $\cW_1$ and $\cW_2$, both of size $q+1$, and the action of $Gal(M|K)\cong \PGL(2,q)$ on $\cW_i$ is the same as on $\Delta$.  From the Hurwitz genus formula \cite[Theorem 11.72]{HKT} applied to the Galois covering $M|K$,
\begin{equation}
\label{eqA291122BKS}
\begin{cases}
2\mathfrak{g}(M)-2&= -2|\PGL(2,q)|+\sum_{T\in \cW}\sum_{i\ge 0}(|G_T^{(i)}|-1)\\
&=-2(q^3-q)+2(q+1)(q(q-1)-1+s(q-1)).
\end{cases}
\end{equation}
The subfield $K(u)$ of $M$ is a Galois cover of $M$ with Galois group $Gal(M|K(u)$ isomorphic to $\AGL(1,q)$. Since $K(u)$ is the function field of $\cC$ which is a rational curve, we have $\mathfrak{g}(K(u))=0$,
the Hurwitz genus formula \cite[Theorem 11.72]{HKT} applied to the Galois covering $M|K(u)$ gives
\begin{equation}
\label{eqB291122BKS}
\begin{cases}
2\mathfrak{g}(M)-2&=-2|\AGL(1,q)|+\sum_{T\in \cW}\sum_{i\ge 0}(|G_T^{(i)}|-1)\\
&=-2q(q-1)+2[(q(q-1)-1+s(q-1)+q(q-2)].
\end{cases}
\end{equation}
Comparison of (\ref{eqA291122BKS}) with (\ref{eqB291122BKS}) yields $s=1$ whence $2\mathfrak{g}(M)-2=2q^2-2q-4$ follows.

Assume now that $m_1=m_2$. Then $\cW_1=\cW_2$, and $\cW_1$ is the only short orbit of $Gal(M|K)$ on $M$. The above computation gives
\begin{equation}
\label{eqA291122BKS1}
\begin{cases}
2\mathfrak{g}(M)-2&= -2|\PGL(2,q)|+\sum_{T\in \cW}\sum_{i\ge 0}(|G_T^{(i)}|-1)\\
&=-2(q^3-q)+(q+1)(q(q-1)-1+s(q-1)).
\end{cases}
\end{equation}
and
\begin{equation}
\label{eqB291122BKS1}
\begin{cases}
2\mathfrak{g}(M)-2&=-2|\AGL(1,q)|+\sum_{T\in \cW}\sum_{i\ge 0}(|G_T^{(i)}|-1)\\
&=-2q(q-1)+(q(q-1)-1+s(q-1)+q(q-2).
\end{cases}
\end{equation}
From this, $s=q+1$ follows. Therefore $2\mathfrak{g}(M)-2=q^2-q-2$.

In the special case, $m_1\ne m_2$ with $m_1\notin \mathbb{F}_q,\, m_2\in \mathbb{F}_q$  and $P_1\ne Q_1, P_2=Q_2$. Let $\cW_2$ be the set of the places of $M$ lying over $P_2$. Since there exists a unique place of $K(u,v)$ lying over $P_2$, $[M:K]=|\PGL(2,q)|=q(q+1)(q-1)$  together with $K(u,v):K]=q+1$ yield $|\cW_2|\le q(q-1)$. Therefore $1$-point stabilizer $G_1$ of $Gal(M:K)\cong \PGL(2,q)$ in $\cW_2$ has order at least $q+1$, say $\lambda(q+1)$ with a divisor $\lambda$ of $q(q-1)$. Actually, $\lambda=1$. In fact, $G_1$ is always solvable and its subgroups of order prime to $p$ are cyclic, see \cite[Lemma 11.44]{HKT}, hence the claim follows from \cite[Theorem A.8]{HKT} which is a corollary to the Dickson's classification of subgroups of $PSL(2,q)$. Therefore $G_1$ is a cyclic group of order $q+1$, and $|\cW_2|=q(q-1)$. Thus $Gal(M|K)$ has two short orbits, namely $\cW_1$ and $\cW_2$ where $|\cW_1|=q+1$ and $Gal(M|K)$ acts on $\cW_1$ as its $3$-transitive permutation representation on $\PG(1,q)$ whereas $\cW_2|=q(q-1)$ and $Gal(M|K)$ acts on $\cW_2$ as on the sets consisting of its cyclic subgroups of order $q+1$, and the action is by conjugacy.
From the Hurwitz genus formula \cite[Theorem 11.72]{HKT} applied to the Galois covering $M|K$,
\begin{equation}
\label{eqA291122BKS2}
\begin{cases}
2\mathfrak{g}(M)-2&= -2|\PGL(2,q)|+\sum_{T\in \cW}\sum_{i\ge 0}(|G_T^{(i)}|-1)\\
&=-2(q^3-q)+(q+1)(q(q-1)-1+s(q-1))+q^2(q-1).
\end{cases}
\end{equation}
and
\begin{equation}
\label{eqB291122BKS2}
\begin{cases}
2\mathfrak{g}(M)-2&=-2|\AGL(1,q)|+\sum_{T\in \cW}\sum_{i\ge 0}(|G_T^{(i)}|-1)\\
&=-2q(q-1)+(q(q-1)-1+s(q-1)+q(q-2).
\end{cases}
\end{equation}
This is only possible for $s=1$. Therefore $2\mathfrak{g}(M)-2=0$, and hence $\mathfrak{g}(M)=0$.

Thus the following claim is proven.
\begin{lemma}
\label{lem291122BKS} In the general case, either
$\mathfrak{g}(M)=q^2-q-1,$ or $\mathfrak{g}(M)=\ha(q^2-q)$ according as $m_1\ne m_2,$ or $m_1=m_2.$  In the special case, $M$ is a rational function field.
\end{lemma}
We are in a position to prove the following theorem.
\begin{theorem}
\label{teoF} Let $k=q^3+1-\ha q(q-1)$. In $\PG(2,q^r)$ with $r\ge 5$, let $\Omega$ be the $(k,q+1)$-arc consisting of all points of the rational BKS curve. For even $r$, $\Omega$ is complete. If $r$ is odd then the points which are uncovered by the $(q+1)$-secants to $\Omega$ are exactly the points in $\PG(2,q)$ not lying in $\Omega$. Adding those points to $\Omega$ produces a complete $(k,q+1)$-arc in $\PG(2,q^r)$, $r$ odd, with $k=q^3+q+1$.
\end{theorem}
\begin{proof} The $\rm{\check{C}}$ebotarev type argument in the proof of Theorem \ref{teoE} in Section \ref{vWvher} can be used to deal with both the general and special cases.

In the general case, if  $m_1\ne m_2$,  the Hasse-Weil lower bound \cite[Theorem 9.18]{HKT} ensures at least $q^{r}+1-2q^{r/2}(q^2-q-1)=q^{r}-2q^{r/2+2}+2q^{r/2+1}+2q^{r/2}$ places of $M$ defined over $\mathbb{F}_{q^{r}}$. For $r>3$, this number is larger than $(q+1)+q(q-1)=q^2+1$ which is the number of ramified places of $M$ under the action of $Gal(M|K)$. Therefore, as long as $r\ge 5$, $M$ has a unramified place $P_0$ over $\mathbb{F}_{2r}$. From this,  as in the proof of Theorem \ref{teoE}, the claim follows. If $m_1=m_2$, the same computation with the Hasse-Weil lower bound proves the existence of a unramified place of $M$ defined over $\mathbb{F}_{q^r}$ as far as $q^r+1-q^{r/2+2}-q^{r/2+1}-(q+1)>0$, i.e. $q\ge 5$.

In the special case, $M$ is rational and hence it has exactly $q^r+1$ places. Furthermore,  $(q+1)+q^2-q=q^2+1$ is the number of the places of $M$ which are ramified under the action of $Gal(G|K)$. Therefore, for $r\ge 3$, $M$ has a unramified place over $\mathbb{F}_{r}$ and the claim can be proven as in the general cases.

We are left with the case where $m_1,m_2\in \mathbb{F}_q$ and $Q\in \PG(2,q)\setminus \Omega$. As mentioned in Section \ref{BKSs}, a linear automorphism group $G\cong \PGL(2,q)$ of $\PG(2,q)$ leaves $\cC$ invariant and acts transitively on the external points to $\cC$ in $\PG(2,q)$. Moreover, the line at infinity meets $\cC$ only in $X_\infty=(1:0:0)$.
Therefore, we may assume that $Q$ is the point at infinity $Y_\infty=(0:1:0)$ so that the line at infinity is not a $(q+1)$-secant to $\Omega$. We show that $Y_\infty$ is covered by a $(q+1)$-secant to $\Omega$ in $\PG(2,q^r)$ if and only if $r$ is even. Let $P\in \PG(2,q^r)$ be a point of $\cC$ with parameter $t_1\in \mathbb{F}_{q^r}$. Then the vertical line through $P$ meets $\cC$ in the points $P_i$ with parameters $t_i$ such that $(t_i+1)^{q+1}=(\tau+1)^{q+1}$. Then $t+1=\lambda (\tau+1)$ with $\lambda^{q+1}=1$. Therefore, $P_i\in \PG(2,q^r)$ for $i=1,2,\ldots, q+1$ if and only if $q+1$ divides $q^r-1$ whence the claim follows.
\end{proof}

%where $G_P{(i)}$ is the $i$-th ramification group of $G$ at $P$. From Serre's result \cite[Lemma 11.75(i)]{HKT}, $G_P^{(0)}$ =$ G_P^{(1)}=\ldots=G_P^{(q-2)}$ is the subgroup of order $q$ and $G_P^{(q-1}=\{1\}$.
\section{Case $P\in \Omega$}
\label{ssin}
We keep our notation from Sections \ref{ssHon} and \ref{ssBKSon}.
As in Section \ref{ssnotin}, our first step is to show if the irreducibility condition in Theorems \ref{teoB} and \ref{teoD} is fulfilled. For this purpose, it is sufficient the $2$-transitivity of $Gal(g,K)$ one the set $\Delta=\{\alpha_1,\ldots,\alpha_q\}$ of the roots of the polynomial $g(T)$. Since $Gal(g,K)$ is transitive by Lemmas \ref{lemA291122} and \ref{lem11122022A}, we may limit ourselves to investigate whether the $1$-point stabilizer of $Gal(g,K)$ is transitive.
\subsection{Hermitian case}
We may assume that $P\notin \PG(2,q^2)$. From Section \ref{hhc}, $P$ is incident with two tangent lines to $\cH_q$, say $\ell_{m_1}$ and $\ell_{m_2}$ where $P$ is the tangency point of $\ell_{m_1}$ while $\ell_{m_2}$ is tangent to $\cH_q$ at the point $R$ whose Frobenius image is $P$, i.e. $R=R(\alpha,\beta)$ with $\alpha^{q^2}=a, \beta^{q^2}=b$. Moreover, $\ell_1$ meets $\cH_q$ in another point $P'$ which is the Frobenius image of $P$. Since $P=P(a,b)$, we have $a^q+m_1=0$ and $a+m_2^q=0$.

Therefore, exactly two places of $K(u)$ lie over $m_1$. They are identified by $P_1=(m_1,0)$ and $P_2=(m_1,-(a+m_1^q))$ where $e(P_{1}|m_1)=q-1$ and $e(P_{2}|m_1)=1$. Instead, there exists just one place of $K(u)$ lying over $m_2$ identified by $R_1=(m_2,-\sqrt[q]{a^q+m_2})$ where $e(R_{1}|m_2)=q$. Furthermore, there are $(q-1)$ places of $M$ over $P_1$ identified by $T_{1,i}=(m_1,0, v_1^i)$ where $v_1$ is a primitive ($q-1$)-st root of
unity. Also, there exists a unique place of $M$ lying over $P_2$ identified by $T_2=(m_1,-(a+m_1^q),0)$, and unique place over $R$ identified by $T_3=(m_2,-\sqrt[q]{a^q+m_2},\infty)$.

From Result \ref{vWsatz}, the inertia group $I(T_2|m_1)$ has two orbits on the set $\Delta$ of the roots of $g(T)$, one of size $q-1$ and another of size $1$. Therefore, $I(T_1|m_1)$ as a subgroup of $Gal(g,K)$ acting on $\Delta$ fixes a point and transitive on the remaining points. Thus $Gal(g,K)$ is $2$-transitive on $\Delta$. Therefore $g_1(T)$ is irreducible over $K$ and $Gal(M|K)$ acts on $\Delta$ as a sharply $2$-transitive permutation group (of order $q(q-1)$).

Let $\cW_1=\{T_{1,i},T_2|i=0,1\ldots, q-2\}$ and $\cW_2=\{T_3\}$. Then $\cW_1$ and $\cW_2$ are the only short orbits of $Gal(M|K)$. Since $Gal(M|K)\cong \AGL(1,q)$, $Gal(M|K)$ has a cyclic subgroup $\Lambda$ of order $q-1$ that fixes a root of $g(T)$. We may assume that $u$ is that root. Then the fixed field of $\Lambda$ is $K(u)=\mathbb{F}(m,u)$. As we have shown, $\mathbb{F}(m,u)$ has a non-singular plane model of degree $q$, and hence $\mathfrak{g}(K(u))=\ha q(q-1)$. Furthermore, $\Lambda$ has exactly two fixed places, namely $T_2$ and $T_3$. From the Hurwitz genus formula \cite[Theorem 11.72]{HKT} applied to the Galois covering $M|K(u)$,
\begin{equation}
\label{eqC291122}
2\mathfrak{g}(M)-2= (q-1)(q(q-1-2)+2(q-2)=q(q-1)^2-2.
\end{equation}
Now, if  $$q^{2r}+1> 2\mathfrak{g} q^r+q+1>q^{r+3}-2q^{r+2}+q^{r+1}+q+1,$$ then the Hasse-Weil lower bound yields that $P$ lies on a $(q+1)$-secant to $\Omega$. In fact, the arguments in the proof of Theorem \ref{teoE} also work for this case and provide a proof for the following result.
\begin{theorem}
\label{teoE1} Let $k=q^{2r}+1\pm q^{r+1}(q-1)$ where $\pm$ is taken according as $r$ is even or odd. In $\PG(2,q^{2r})$ with $r\ge 3$, let $\Omega$ be the $(k,q+1)$-arc consisting of all points of the Hermitian curve. If $r>3$ then each point of $\Omega$ is covered by some $(q+1)$-secant to $\Omega$.
\end{theorem}
For $r=3$, the results stated in Section  \ref{ssHIII} are sufficient to determine the number $n_P$ of ($q+1$)-secants to $\cH_q$ through any point $P\in PG(2,q^6)\setminus \PG(2,q^2)$ lying in $\cH_q$. We may assume $\cH_q$ in its canonical form $X^q+T^qX+T=0$ over $\mathbb{F}_{q^6}$. Since the automorphism group of $\cH_q$ has two orbits in $PG(2,q^6)$, one consisting of its points in $PG(2,q^2)$, we may also assume $P=X_\infty$. As the Galois closure $M$ is a maximal curve of genus $\ha q(q-1)^2$, it has as many as $q^6+1+q(q-1)^2q^3$ points in $PG(2,q^6)$. Moreover, the Galois extension $M|L(u)$ has degree $q(q-1)$ where $L(u)$ is the function field of $\cH_q$.  Therefore, $n_Pq(q-1)=q^6+1+q(q-1)^2q^3-(q+1)^2$ whence $n_P=2q^4+q^2+q+1$ follows. A direct proof based on norms and traces in finite fields is also possible. This was pointed out by B. Csajb\'ok \cite{Cs}.

\subsection{The rational BKS case}
We begin with case $P\not\in \PG(2,q)$. From Section \ref{BKSs}, $P$ is incident with two tangent lines to $\cC$, say $\ell_{m_1}$ and $\ell_{m_2}$ where $P$ is the tangency point of $\ell_{m_1}$ while $\ell_{m_2}$ is tangent to $\cC$ at the point $R$ whose Frobenius image is $P$. As in Section \ref{ssHon}, we replace $2m/(1-2m)$ by $m$ so that we may use equation (\ref{eq3jul2022C1}). Then, if $P$ has parameter $t\in \mathbb{F}$, then the tangent line $\ell_1$ at $P$ has slope $-t^q$ and meets $\cC$ in another point, namely that with parameter $t^q$.

Thus, exactly two places of $K(u)$ lie over $t$. They are identified by $P_1=(-t^q,0)$ and $P_2=(-t^q,t^q-t)$. Moreover, there exists just one place over $-t$ identified by $R_1=(-t,\sqrt[q]{t-t^q})$. Furthermore, there are $(q-1)$ places of $M$ over $P_1$ identified by $T_{1,i}=(-t^q,0,-v_1^i)$ where $v_1$ is a primitive ($q-1$)-st root of
unity. Also, there exists a unique place of $M$ lying over $P_2$ identified by $T_2=(-t^q,t^q-t,0)$, and unique place over $R$ identified by $T_3=(-t,-\sqrt[q]{t-t^q},\infty)$.
This ramification picture is analog of that in Section \ref{ssin}. In particular, Result \ref{vWsatz} yields that the inertia group $I(T_2|-t)$ has two orbits on the set $\Delta$ of the roots of $g(T)$, one of size $q-1$ and another of size $1$. Therefore, $I(T_1|-t)$ as a subgroup of $Gal(g,K)$ acting on $\Delta$ fixes a point and transitive on the remaining points. Thus $Gal(M|K)$ is $2$-transitive and $g_1(T)$ is irreducible over $K$.
More precisely, $Gal(g,K)$ acts on $\Delta$ as a sharply $2$-transitive permutation group (of order $q(q-1)$). As in Section \ref{ssin}, we introduce $\cW_1,\cW_2,\Lambda$ and $u$. This time the fixed field of $\Lambda$, that is, $\mathbb{F}(m,u)$ given in (\ref{eqA11122022Y}) which a rational curve.  From the Hurwitz genus formula \cite[Theorem 11.72]{HKT} applied to the Galois covering $M|K(u)$,
\begin{equation}
\label{eqC291122BKS}
2\mathfrak{g}(M)-2=-2(q+1)+2(q-2)=-2.
\end{equation}
Thus $M$ is a rational function field, and hence it has exactly $q^r+1$ places over $\mathbb{F}_q$. Therefore, the number of long orbits of $Gal(M|K)$ defined over $\mathbb{F}_{q^r}$ is equal to $(q^r+1-(q+1))/(q(q-1))=1+q+\ldots+q^{r-2}$.

Now, suppose that $P\in \PG(2,q)$. If $P$ is an internal point to $\cC^2$ then $P$ is a singular point. Clearly, no singular point of $\cC$ is covered by a $(q+1)$-secant to $\cC$ as $\deg(\cC)=q+1$. Therefore,
we are left with the case where $P\in \cC^2$. As mentioned in Section \ref{BKSs}, a linear automorphism group $G\cong \PGL(2,q)$ of $\PG(2,q)$ leaves $\cC$ invariant and acts transitively on the points of $\cC^2$ in $\PG(2,q)$
Therefore, we may assume $P=X_\infty(1:0:0)$. Take point $R\in \cC$ with parameter $t_0\in \mathbb{F}_{q^r}\setminus \mathbb{F}_{q^{2}}$. The horizontal line $\ell$ through $R$ meets $\cC$ in the points with parameters $t$ such that $t^q+t+2=t_0^q+t_0=2$, that is, $(t-t_0)^q+(t-t_0)=0$. Then $t=t_0+\tau$ with $\tau\in \mathbb{F}_q$. Hence $t\in \mathbb{F}_{q^r}\setminus \mathbb{F}_{q^{2}}$. Thus the common points of $\ell$ and $\cC$ are pairwise distinct and each lies in $\PG(2,q^r)$.

Therefore the following result is obtained.
\begin{theorem}
\label{teoF1} Let $k=q^3+1-\ha q(q-1)$. In $\PG(2,q^r)$ with $r\ge 3$, let $\Omega$ be the $(k,q+1)$-arc consisting of all points of the rational BKS curve. Then the points of $\Omega$ uncovered by the $(q+1)$-secants to $\Omega$ are exactly the points of $\PG(2,q)$ lying in $\Omega$.
\end{theorem}

\section{Acknowledgement} The research was supported by the Hungarian Academy of Science ``MTA Vendégkutatói Program 2022'', the NKFIH-OTKA Grants SNN 132625, K 124950, and the Program of Excellence TKP2021-NVA-02 at the Budapest University of Technology and Economics. The third author was partially supported by the Slovenian Research Agency, research project J1-9110.

\end{document}